\documentclass[1p,times,number]{elsarticle}
\usepackage{amsmath,amssymb,amsfonts,amsthm}
\usepackage{graphicx}
\usepackage{booktabs}
\usepackage{xurl}
\usepackage[T1]{fontenc}
\usepackage{microtype}
\biboptions{sort&compress}
\newtheorem{theorem}{Theorem}[section]
\newtheorem{lemma}[theorem]{Lemma}
\newtheorem{proposition}[theorem]{Proposition}
\newtheorem{corollary}[theorem]{Corollary}
\theoremstyle{definition}
\newtheorem{example}[theorem]{Example}

\theoremstyle{remark}
\newtheorem{remark}[theorem]{Remark}
\numberwithin{equation}{section}

\journal{Journal of Algebra}
\begin{document}
\begin{frontmatter}
\title{Iso-Multiplication Modules over Noncommutative Rings:\\
Triangular Classification and Rational Rigidity}
\author[mpei]{Askar Tuganbaev\corref{cor1}}
\ead{tuganbaev@gmail.com}
\author[marmara]{Unsal Tekir}
\ead{utekir@marmara.edu.tr}
\author[marmara]{Suat Koc}
\ead{suat.koc@marmara.edu.tr}
\author[gebze]{Gulsen Ulucak}
\ead{gulsenulucak@gtu.edu.tr}
\cortext[cor1]{Corresponding author.}
\address[mpei]{Department of Higher Mathematics, National Research University MPEI, Moscow, Russia}
\address[marmara]{Department of Mathematics, Marmara University, Istanbul, T\"urkiye}
\address[gebze]{Department of Mathematics, Gebze Technical University, Kocaeli, T\"urkiye}
\begin{abstract}
A right $A$-module $M$ is called iso-multiplication if for every submodule
$N$ of $M$ there is a two-sided ideal $B$ of $A$ such that $N\cong_A MB$.
We construct noncommutative examples of iso-multiplication modules which
are not multiplication modules, including triangular algebras of arbitrary
dimension. We classify all iso-multiplication modules over $T_2(k)$ without
a finite-generation assumption. For a semiprime right Goldie ring we prove
that the classical localization of an iso-multiplication module is
multiplicity-free. We also obtain restrictions on the composition length,
results for local rings, and connections with Morita equivalence and chain
conditions up to isomorphism.
\end{abstract}
\end{frontmatter}

\section{Introduction}

All rings are assumed to be associative and to have an identity $1\neq0$;
all modules are unital right modules. An ideal means a two-sided ideal.
For a right $A$-module $M$ and an ideal $B$ of $A$, we denote by $MB$ the
submodule consisting of all finite sums of elements $mb$, where $m\in M$
and $b\in B$. Following Tuganbaev \cite{Tuganbaev2003,Tuganbaev2004}, the
module $M$ is called a \emph{multiplication module} (an m-module) if for every submodule
$N$ of $M$ there exists an ideal $B$ of $A$ such that $N=MB$.

We call $M$ an \emph{iso-multiplication module} (an im-module) if for every submodule
$N\le M$ there exists an ideal $B$ of $A$ such that $N\cong_A MB$.
The isomorphism is $A$-linear and need not be the restriction of an
automorphism of $M$. Every multiplication module is an iso-multiplication module.
For noncommutative rings the converse is not true. In particular, a right
ideal can be isomorphic to a two-sided ideal without being a two-sided ideal.

The principal examples are the regular right modules of triangular rings.
For a field $k$, the right regular module of $T_2(k)$ is iso-multiplication
but is not a multiplication module. We give the complete classification of
iso-multiplication modules over this ring, including infinitely generated
modules. The triangular construction extends to the rings $A(V)$ associated
with arbitrary nonzero $k$-vector spaces $V$; it also gives examples of
infinite dimension. Related examples are obtained from skew polynomial
rings, direct products, factor rings, and Morita equivalences.

We prove that an iso-multiplication module over a semiprime right Goldie
ring has multiplicity-free classical localization. In particular, if the
regular module $A_A$ is iso-multiplication, the classical right ring of
quotients of $A$ is a finite product of division rings. This condition is
necessary but not sufficient. We also obtain bounds on composition length
when the ring has finitely many ideals, a cyclicity result for local rings,
and results on chain conditions up to isomorphism.

We give some notation. The set of submodules of $M$ is denoted by
$\mathcal{L}(M)$. For $N\le M$ we set
\[
(N:M)=\{a\in A:Ma\subseteq N\}=r_A(M/N).
\]
Thus $(N:M)$ is a two-sided ideal. We denote by $\operatorname{Ann}_A(X)$
the annihilator of a right module $X$, by $l_A(X)$ its composition length
when finite, and by $\operatorname{u.dim}_A(X)$ its right uniform dimension.
For a semiprime right Goldie ring $A$, we write $C$ for the set of regular
elements, $Q_{\mathrm{cl}}^r(A)$ for its classical right ring of quotients,
and $M_Q=M\otimes_AQ$ when $Q=Q_{\mathrm{cl}}^r(A)$. The terms
\emph{isonoetherian} and \emph{isoartinian} are used for the corresponding
chain conditions up to isomorphism \cite{Facchini2016,Facchini2017}.

\section{Examples and elementary properties}

For a right $A$-module $M$ and a two-sided ideal $B$ of $A$, we use
$MB$ in the sense defined in the Introduction. In particular, for the
regular module $A_A$, we have $AB=B$. The distinction between right
ideals and two-sided ideals is essential in the following examples.

\begin{proposition}
\label{prop:precisely_changes} For every right module
$M$ the following hold:

\begin{enumerate}
\item[(i)] $M$ is a multiplication module if and only if $N = M(N:M)$ for every $N
\in\mathcal{L}(M)$.

\item[(ii)] Every multiplication module is an iso-multiplication module; the iso-multiplication property is invariant under isomorphisms.

\item[(iii)] The zero module and every simple module are iso-multiplication modules (indeed, multiplication modules).
\end{enumerate}
\end{proposition}

\begin{proof}
If $N = MB$, then $B \subseteq(N:M)$, so $N = MB \subseteq M(N:M) \subseteq
N$. Conversely, use $B = (N:M)$. The other assertions follow from the
definitions; an isomorphism $f: M \to M^{\prime}$ maps $MB$ onto $M^{\prime}B$
while a simple module has only the submodules $0 = M0$ and $M = MA$.
\end{proof}

\begin{remark}
For an iso-multiplication module, $N \cong MB$ does not imply $B \subseteq(N:M)$, and it does
not imply $N \cong M(N:M)$. The distinction is demonstrated explicitly in
Example \ref{ex:basic_counterexample}.
\end{remark}

\begin{proposition}
The module $M_{A}$ is an iso-multiplication module precisely
when the map
\[
\{B : B \triangleleft A\} \longrightarrow\{[N] : N \in\mathcal{L}(M)\}, \quad
B \mapsto[MB]
\]
is surjective. In particular, for every $N \le M$ there is an ideal $B$ with
$\text{Ann}_{A}(N) = \text{Ann}_{A}(MB)$; equality of annihilators alone is
not sufficient.
\end{proposition}

\begin{proof}
The first assertion is the definition expressed in terms of isomorphism
classes. An $A$-linear isomorphism preserves annihilators. To see why annihilator equality alone does not imply the defining
isomorphism, let $A=k$ and $M=k^{2}$. The submodule $N=k\oplus0$ and the
ideal multiple $MA=M$ have the same zero annihilator, but
$N\not\cong_A M$.
\end{proof}

\subsection{Examples over triangular and principal right ideal rings}

\begin{example}
\label{ex:basic_counterexample} Let $k$ be any field, $A = T_{2}(k) = \left\{
\begin{pmatrix}
a & b\\
0 & c
\end{pmatrix}
: a,b,c \in k \right\}  $ and $M = A_{A}$. Denote the matrix units by $e_{11},
e_{12}, e_{22}$. The complete list of two-sided ideals is
\[
0, \quad J = ke_{12}, \quad K = ke_{12} \oplus ke_{22}, \quad L = ke_{11}
\oplus ke_{12}, \quad A.
\]
The right $A$-module $K$ is a direct sum of two copies of $S_{2} = e_{22}A =
ke_{22}$; both $e_{12}$ and $e_{22}$ are multiplied by $%
\begin{pmatrix}
a & b\\
0 & c
\end{pmatrix}
$ through the scalar $c$. Consequently, any right ideal $U \subseteq K$ is a
$k$-linear subspace and satisfies
\[
U \cong%
\begin{cases}
0, & \dim_{k} U = 0,\\
J, & \dim_{k} U = 1,\\
K, & \dim_{k} U = 2.
\end{cases}
\]
If a right ideal $U$ is not contained in $K$, choose $u = ae_{11} + be_{12} +
ce_{22} \in U$ with $a \neq0$. Then $ue_{11} = ae_{11}$ and $e_{11}e_{12} =
e_{12}$ show $L = e_{11}A \subseteq U$. Since $A/L$ is simple, $U = L$ or $A$.
Every right ideal is therefore isomorphic to one of $0, J, K, L, A$. Since $MB
= AB = B$ for each ideal $B$, $A_{A}$ is iso-multiplication.

It is not a multiplication module: $N = e_{22}A = ke_{22}$ is a right ideal but not
two-sided, because $e_{12}e_{22} = e_{12} \notin N$. Nevertheless,
\[
e_{22}A \longrightarrow J = e_{12}A, \quad e_{22}a \mapsto e_{12}a
\]
is a well-defined $A$-isomorphism. Moreover $(N:M) = 0$, as there is no
nonzero two-sided ideal inside $ke_{22}$. Thus $M(N:M) = 0 \neq N$.
\end{example}

\begin{table}[htbp]
\caption{Isomorphism types of right ideals of $T_{2}(k)$}%
\centering
\vspace{2mm}
\begin{tabular}
[c]{ll}%
\toprule Isomorphism type of a right ideal of $T_{2}(k)$ & A representing
two-sided ideal\\
\midrule $0$ & $0$\\
$S_{2}$ (any line in $K$) & $J = ke_{12}$\\
$S_{2} \oplus S_{2}$ & $K = ke_{12} \oplus ke_{22}$\\
$P = e_{11}A$ & $L = e_{11}A$\\
$A_{A} = P \oplus S_{2}$ & $A$\\
\bottomrule
\end{tabular}
\end{table}

\begin{proposition}
\label{prop:domain_class} If $A$ is a (possibly
noncommutative) domain such that every right ideal is principal, then $A_{A}$
is iso-multiplication. It is a multiplication module if and only if every right ideal is two-sided.
\end{proposition}

\begin{proof}
For every nonzero right ideal $aA$, multiplication $A_{A} \to aA$, $r \mapsto
ar$, is an isomorphism: surjectivity is clear and injectivity follows from the
absence of zero divisors. Thus each right ideal is isomorphic to $A = AA$ or
$0 = A0$. For the second assertion use $AB = B$ for a two-sided ideal $B$.
\end{proof}

\begin{example}
\label{ex:skew_polynomial} Let
$K = k(t)$ with $\text{char } k = 0$, let $\sigma(t) = t+1$ on $K$, and
consider the Ore extension $A = K[x;\sigma]$ subject to $xa = \sigma(a)x$. Its
degree function is additive on nonzero products; in particular $A$ is a
domain. Here every nonzero right ideal is principal. Indeed, let $I$ be such an
ideal and choose $0\ne g\in I$ of minimal degree $d$. Given $f\in I$
with $\deg f\ge d$, write $m=\deg f-d$. If $a$ and $b$ are the leading
coefficients of $g$ and $f$, respectively, choose
$c=\sigma^{-d}(a^{-1}b)\in K$. The leading term of $g(cx^{m})$ is
$ax^{d}cx^{m}=a\sigma^{d}(c)x^{d+m}=bx^{\deg f}$. Subtracting
$g(cx^{m})$ from $f$ therefore lowers the degree. Repeating the step gives
$f=gq+r$ with $\deg r<d$ (or $r=0$). Since $r\in I$, minimality of $d$
forces $r=0$. Hence $I=gA$, so $A_A$ is iso-multiplication by Proposition
\ref{prop:domain_class}.

The right ideal $(x-1)A$ is not two-sided. Indeed, if it were a left ideal,
then $t(x-1)=(x-1)f$ for some $f\in A$. Since $A$ is a domain with additive
degree on nonzero products, $f$ must have degree zero, say $f=c\in K$.
Expanding with $xc=\sigma(c)x$ gives
\[
tx-t=(x-1)c=\sigma(c)x-c.
\]
The constant coefficient on the left is $-t$ and on the right
is $-c$, so $c=t$. The coefficient of $x$ on the left is $t$ and on
the right is $\sigma(c)$, so $\sigma(c)=t$. Substitution gives
$t+1=\sigma(t)=t$, a contradiction. Hence $A_{A}$ is not a classical multiplication module. This
supplies an example over a non-Artinian, noncommutative domain.
\end{example}

\begin{example}
For every $r \ge0$, let $A_{r} =
T_{2}(k) \times k^{r}$. Then the right regular module $(A_{r})_{A_{r}}$ is iso-multiplication
but not a multiplication module. Indeed, the $T_{2}(k)$ factor is iso-multiplication by Example
\ref{ex:basic_counterexample}; each regular $k$-factor is simple, and the
finite-product theorem proved below (Theorem \ref{thm:finite_direct_products})
applies. The non-multiplication property is inherited from the factor
$T_{2}(k)$ by considering the right ideal $e_{22}T_{2}(k) \times0$.
\end{example}

\begin{example}
\label{ex:submodule_failure} In Example \ref{ex:basic_counterexample}, the
right ideal $K = ke_{12} \oplus ke_{22}$ is a submodule of the iso-multiplication module
$A_{A}$, but $K_{A}$ is not iso-multiplication. Indeed $K \cong S_{2} \oplus S_{2}$. The only
ideal multiples of $K$ are $0$ and $K$: the action of $A$ on $K$ factors
through the bottom-right diagonal entry, and each two-sided ideal either has
zero bottom-right entry or contains $e_{22}$. The proper nonzero submodule $J
\cong S_{2}$ is isomorphic to neither $0$ nor $K$.
\end{example}

\begin{example}
\label{ex:direct_sum_failure}
For any division ring $D$, the simple right module $D_{D}$ is a multiplication module, but
$D_{D} \oplus D_{D}$ is not iso-multiplication: the only ideal multiples are $0$ and $D^{2}$,
whereas a copy of $D$ is a submodule. In particular the iso-multiplication property is not
closed under finite direct sums.
\end{example}

\begin{remark}
The iso-multiplication condition is not automatic, even in the commutative
setting: for the Pr\"ufer group $\mathbb{Z}(p^{\infty})$, every nonzero ideal multiple equals the
whole group, whereas it contains a submodule of order $p$. We do not pursue
commutative classifications here.
\end{remark}

\subsection{Triangular rings of arbitrary dimension}

The preceding construction extends to triangular rings of arbitrary
dimension. We describe their right ideals explicitly.

\begin{theorem}
\label{thm:triangular_family} Let $k$ be a
field, $0 \neq V$ an arbitrary $k$-vector space, and put
\[
A(V) = \left\{
\begin{pmatrix}
a & v\\
0 & c
\end{pmatrix}
: a,c \in k, v \in V \right\}  ,
\]
\[%
\begin{pmatrix}
a & v\\
0 & c
\end{pmatrix}
\begin{pmatrix}
a^{\prime} & v^{\prime}\\
0 & c^{\prime}%
\end{pmatrix}
=
\begin{pmatrix}
aa^{\prime} & av^{\prime}+ vc^{\prime}\\
0 & cc^{\prime}%
\end{pmatrix}
.
\]
Here $V$ is regarded as the natural $(k,k)$-bimodule, with the usual scalar
multiplication on both sides; this convention is essential for the
classification of two-sided ideals. Then the regular right module
$A(V)_{A(V)}$ is iso-multiplication but not a classical multiplication module.
More precisely, writing
\[
J = \left\{
\begin{pmatrix}
0 & v\\
0 & 0
\end{pmatrix}
: v \in V \right\}  , \quad K = J \oplus ke_{22}, \quad L = e_{11}A(V),
\]
every right ideal is either a $k$-subspace of $K$, or $L$, or $A(V)$; each
subspace of $K$ is isomorphic to a two-sided ideal of $A(V)$.
\end{theorem}

\begin{proof}
The action of $A(V)$ on $K$ factors through the homomorphism $A(V) \to k$
taking a matrix to its bottom-right entry. Hence every $k$-subspace $U
\subseteq K$ is a right ideal, and $U \cong S_{2}^{(\dim_{k} U)}$, where
$S_{2} = ke_{22}$ and the parenthesized exponent permits infinite cardinals.
Since $V \neq0$, for each cardinal $\lambda\le\dim_{k} K$ there exists a
two-sided ideal isomorphic to $S_{2}^{(\lambda)}$: if $\lambda\le\dim_{k} V$,
choose a $\lambda$-dimensional subspace $W \subseteq V$ and use $W \subseteq
J$ which is a two-sided ideal; otherwise necessarily $\dim_{k} V$ is finite
and $\lambda= \dim_{k} V + 1$, in which case use the two-sided ideal $K$. (For
infinite $V$, $\dim_{k} K = \dim_{k} V$.)

If a right ideal $U$ is not contained in $K$, it contains a matrix $u$ with
nonzero top-left entry $a$. Then $ue_{11} = ae_{11}$ implies $e_{11} \in U$
and hence $L = e_{11}A(V) \subseteq U$. Since $A(V)/L \cong S_{2}$ is simple,
$U = L$ or $A(V)$. Both $L$ and $A(V)$ are themselves two-sided ideals. Thus
every right ideal is isomorphic to a two-sided ideal, which proves that
$A(V)_{A(V)}$ is iso-multiplication.

Finally $e_{22}A(V) = ke_{22}$ is a right ideal but not a two-sided ideal: for
$0 \neq v \in V$ the upper-right matrix $j_{v}$ satisfies $j_{v} e_{22} =
j_{v} \notin ke_{22}$. The regular module is therefore not multiplication.
\end{proof}

\begin{corollary}
\label{cor:iso_chain_separation} If $V$ is infinite dimensional, then
$A(V)_{A(V)}$ is iso-multiplication and isoartinian, but it is not
isonoetherian. In particular, iso-multiplication plus isoartinianity does not
imply isonoetherianity, even for a noncommutative right regular module.
\end{corollary}

\begin{proof}
By Theorem \ref{thm:triangular_family} every right ideal is a subspace of $K$,
or $L$, or $A(V)$. In a descending chain, terms outside $K$ are confined to
the two possibilities $A(V)$ and $L$, so after at most two such terms (up to
repetition) the chain either remains equal to $L$ or enters $K$ permanently.
Once inside $K$, it is a descending chain of vector spaces. Its sequence of
dimension cardinals is nonincreasing and must eventually be constant, since
there is no infinite strictly descending sequence of cardinals. The
corresponding right ideals are then mutually isomorphic as right modules. Thus
the regular module is isoartinian.

For the ascending chain, choose linearly independent vectors
$v_{1},v_{2},\ldots\in V$ and define, for $n\geq 1$,
\[
U_{n}=kv_{1}\oplus\cdots\oplus kv_{n}\subseteq J.
\]
Then $U_{1}\subsetneq U_{2}\subsetneq\cdots$ is a strictly ascending
chain of right ideals (each $U_n$ is a subspace of the two-sided ideal
$J$). The action on $J$ is through the lower-right
scalar, so $U_{n}\cong S_{2}^{\oplus n}$ for every $n$. In particular,
$U_{n}\not\cong U_{m}$ when $n\neq m$, as their composition lengths
are different. This chain never stabilizes up to isomorphism, and hence
$A(V)_{A(V)}$ is not isonoetherian.
\end{proof}

\begin{remark}
For finite-dimensional $V$, the ring $A(V)$ is right Artinian and consequently
both isoartinian and isonoetherian. For infinite-dimensional $V$ it has a
square-zero Jacobson radical $J$, so it lies outside the semiprime right
Goldie setting of Section \ref{sec:rational_rigidity}. No comparison with a
semiprime classification is implied.
\end{remark}

\section{Modules over the triangular algebra $T_2(k)$}

The right modules over $T_2(k)$ are representations of the quiver with one
arrow. We recall their standard decomposition and determine which of these
modules are iso-multiplication, without a finiteness assumption.

Let $A = T_{2}(k)$ as in Example \ref{ex:basic_counterexample}. Put
\[
P = e_{11}A, \quad S_{1} = e_{11}A / e_{12}A, \quad S_{2} = e_{22}A.
\]
Here $P$ is a two-dimensional indecomposable projective, with $\text{rad } P
\cong S_{2}$ and top $S_{1}$. The two simple right modules are $S_{1}$ and
$S_{2}$.

\begin{lemma}
\label{lem:normal_form} Every
right $A$-module, with no finiteness condition, decomposes uniquely up to
isomorphism as
\[
M \cong P^{(a)} \oplus S_{1}^{(b)} \oplus S_{2}^{(c)}, \quad a, b, c \text{
cardinals.}
\]

\end{lemma}

\begin{proof}
Put $V_{1} = Me_{11}$ and $V_{2} = Me_{22}$. Right multiplication by $e_{12}$
gives a $k$-linear map $f: V_{1} \to V_{2}$, $v \mapsto ve_{12}$. Conversely
an arbitrary right $A$-module is determined by this map, because a matrix $%
\begin{pmatrix}
u & w\\
0 & z
\end{pmatrix}
$ acts on $v_{1} + v_{2}$ as $uv_{1} + wv_{1} e_{12} + zv_{2}$. Choose
complements $V_{1} = \text{ker } f \oplus U$ and $V_{2} = f(U) \oplus W$; they
exist even for infinite-dimensional vector spaces. The restriction $f|_{U}: U
\to f(U)$ is an isomorphism. Choose bases of $U$, $\text{ker } f$ and $W$. The
cardinal $a = \dim_{k} U$ gives $P^{(a)}$, while $b = \dim_{k} \text{ker } f$
and $c = \dim_{k} W$ give $S_{1}^{(b)}$ and $S_{2}^{(c)}$. These cardinals are
uniquely determined by the linear map $f$.
\end{proof}

\begin{lemma}
\label{lem:ideal_multiples_T2} For $M = P^{(a)} \oplus
S_{1}^{(b)} \oplus S_{2}^{(c)}$, the ideal multiples of $M$, up to $A$-module isomorphism, are
\[
\begin{array}{c|ccccc}
B&0&J&K&L&A\\ \hline
[MB]&[0]&[S_{2}^{(a)}]&[S_{2}^{(a+c)}]&[P^{(a)}\oplus S_{1}^{(b)}]&[M]
\end{array}
\]
\end{lemma}

\begin{proof}
For the representation $f: V_{1} \to V_{2}$ of Lemma \ref{lem:normal_form},
identify $V_1=Me_{11}$ and $V_2=Me_{22}$ as additive subspaces of $M$.
Then $MJ = \operatorname{im} f$, $MK = V_{2}$, and
$ML = V_{1} \oplus \operatorname{im} f$ as additive subspaces of $M$;
the latter is the subrepresentation consisting of all of $V_{1}$ and
only the image of $f$ in $V_{2}$. These are isomorphic respectively to $S_{2}^{(a)}$, $S_{2}^{(a+c)}$,
and $P^{(a)} \oplus S_{1}^{(b)}$. Here brackets in the table denote
isomorphism classes, not equality of submodules inside $M$; for example
$MJ$ is the image of $f$ inside $Me_{22}$. The cases $0, A$ are immediate.
\end{proof}

\begin{theorem}
\label{thm:T2_classification} Let $A =
T_{2}(k)$ and let $M$ be an arbitrary right $A$-module, possibly infinitely
generated. Then $M$ is iso-multiplication if and only if
\[
M \cong0, \quad S_{1}, \quad S_{2}, \quad P, \quad S_{1} \oplus S_{2},
\quad\text{or} \quad P \oplus S_{2}.
\]
Among these six isomorphism types the first five are multiplication modules, and $P \oplus
S_{2} \cong A_{A}$ is iso-multiplication but not a multiplication module. Thus $A_{A}$ is the only isomorphism type
over $T_{2}(k)$ separating the two conditions. In particular, every iso-multiplication module
over $T_{2}(k)$ is finite dimensional over $k$.
\end{theorem}

\begin{proof}
By Lemma \ref{lem:normal_form}, write $M = P^{(a)} \oplus S_{1}^{(b)} \oplus
S_{2}^{(c)}$, where $a, b, c$ are arbitrary cardinals. Assume first that $M$
is iso-multiplication.

If $b > 0$, the summand $S_{1}$ is a submodule. By Lemma
\ref{lem:ideal_multiples_T2}, the only ideal multiples with any $S_{1}$
composition factor are $ML = P^{(a)} \oplus S_{1}^{(b)}$ and $M$ itself. Thus,
in order that an ideal multiple be isomorphic to $S_{1}$, necessarily $a = 0$
and $b = 1$. Therefore $M = S_{1} \oplus S_{2}^{(c)}$. If $c > 0$, a copy of
$S_{2}$ must also be isomorphic to an ideal multiple; now $MK = S_{2}^{(c)}$
is the only nonzero ideal multiple of pure $S_{2}$-type, so $c = 1$. Hence
this case gives $S_{1}$ or $S_{1} \oplus S_{2}$.

If $b = 0$ and $a = 0$ then $M = S_{2}^{(c)}$. Since all nonzero ideal
multiples equal $M$, the presence of a simple submodule forces $c = 1$ unless
$M = 0$.

Finally suppose $b = 0$ and $a > 0$. The module $M$ contains $S_{2}
\subseteq\text{rad } P$. The only nonzero pure $S_{2}$ ideal multiples in
Lemma \ref{lem:ideal_multiples_T2} are $S_{2}^{(a)}$ and $S_{2}^{(a+c)}$;
hence $a = 1$ (a simple module cannot be isomorphic to an infinite direct sum
of simples). If $c \ge2$, the semisimple submodule $S_{2} \oplus S_{2}
\subseteq S_{2}^{(1+c)} = MK$ has no isomorphic ideal multiple, since the only
nonzero pure $S_{2}$ multiples have multiplicities $1$ and $1+c \ge3$,
including the possibility that $c$ is infinite. Thus $c \le1$. We obtain $P$
and $P \oplus S_{2}$.

Conversely, $0, S_{1}, S_{2}$ are simple or zero and hence m. The only
submodules of $P$ are $0$, $\text{rad } P = PJ \cong S_{2}$, and $P$, so $P$
is a multiplication module. Because $S_{1}, S_{2}$ are nonisomorphic simple modules, the submodules
of $S_{1} \oplus S_{2}$ are $0, S_{1}, S_{2}, M$, and they are respectively
$M0, ML, MK, MA$; hence this module is a multiplication module. Lastly $P \oplus S_{2} \cong e_{11}A
\oplus e_{22}A = A_{A}$. For completeness, its right ideals fall into the
following isomorphism types: $0, S_{2}, S_{2} \oplus S_{2}, P$, and $A$.
Indeed a right ideal contained in $K = ke_{12} \oplus ke_{22}$ is a
$k$-subspace of $K$ and is semisimple of $S_{2}$-type, whereas a right ideal
not contained in $K$ contains $e_{11}$ and hence is $L = e_{11}A$ or $A$.
These five types are represented by the two-sided ideals $0, J, K, L, A$,
respectively. Therefore $A_{A}$ is iso-multiplication; it is not multiplication because $e_{22}A$ is not
two-sided, as shown in Example \ref{ex:basic_counterexample}.
\end{proof}

\begin{table}[th]
\caption{Classification summary over $T_{2}(k)$}%
\centering
{\small \vspace{2mm}
\begin{tabular}
[c]{@{}p{0.22\textwidth}p{0.47\textwidth}p{0.13\textwidth}@{}}%
\toprule Iso-multiplication type & Representative ideal multiples occurring
inside $M$ & Classical m?\\
\midrule $0$ & $0$ & Yes\\
$S_{1}$ & $0, S_{1}$ & Yes\\
$S_{2}$ & $0, S_{2}$ & Yes\\
$P$ & $0, S_{2}, P$ & Yes\\
$S_{1} \oplus S_{2}$ & $0, S_{1}, S_{2}, M$ & Yes\\
$P \oplus S_{2} \cong A_{A}$ & $0, J, K, L, A$ (see Example
\ref{ex:basic_counterexample}) & No\\
\bottomrule
\end{tabular}
}\end{table}

\begin{remark}
The normal form is the standard decomposition of representations of the
one-arrow quiver of type $A_2$. Theorem \ref{thm:T2_classification}
determines which of these representations are iso-multiplication. Among
them, only the regular module fails to be a multiplication module.
\end{remark}

\begin{corollary}
\label{cor:T2_exception} Let
$A=T_{2}(k)$ and let $M$ be an iso-multiplication right $A$-module. Then $M$
is a classical multiplication module if and only if $M \ncong A_{A}$.
\end{corollary}

\begin{proof}
By Theorem \ref{thm:T2_classification}, exactly six isomorphism classes have
the iso-multiplication property. The first five classes are classical
multiplication modules, while $P\oplus S_{2}\cong A_{A}$ is not a classical
multiplication module, by Example \ref{ex:basic_counterexample}.
\end{proof}

\begin{remark}
The decomposition lemma and classification include infinite-dimensional
representations. The classification shows a posteriori that no
infinite-dimensional module over $T_{2}(k)$ satisfies the iso-multiplication property. This
does not extend automatically to triangular algebras with higher-dimensional
off-diagonal bimodules.
\end{remark}

\section{Composition length and local rings}

The number of two-sided ideals of a ring gives restrictions on the
composition lengths of its iso-multiplication modules.

\begin{proposition}
\label{prop:length_spectrum} Let $M_{A}$ be an
iso-multiplication module of finite composition length $r$. Then
\[
\{0, 1, \dots, r\} \subseteq\{l_{A}(MB) : B \triangleleft A\}.
\]
In particular there are at least $r+1$ distinct ideal multiples $MB$. This
conclusion does not imply that the corresponding ideals can be chosen in a
chain under inclusion.
\end{proposition}

\begin{proof}
Choose a composition series $0 = M_{0} \subsetneq M_{1} \subsetneq
\dots\subsetneq M_{r} = M$. For every $j$ the defining condition supplies a
two-sided ideal $B_{j}$ with $M_{j} \cong MB_{j}$. Hence $l_{A}(MB_{j}) = j$.
Since distinct lengths preclude isomorphism, the multiples $MB_{j}$ are
distinct. Nothing in this argument compares $B_{i}$ with $B_{j}$ by inclusion.
\end{proof}

\begin{theorem}
Let $A$ be any unital ring and
$M_{A}$ an iso-multiplication module. Assume every finitely generated submodule of $M$ has
finite composition length. If the collection $\{[MB] : B \triangleleft A\}$
consists of exactly $u < \infty$ isomorphism types, then $M$ has finite length
and
\[
l_{A}(M) \le u-1.
\]
In particular, $M$ is both Noetherian and Artinian, and thus both
isonoetherian and isoartinian. Neither $A$ nor $M$ is assumed Noetherian or
Artinian in advance.
\end{theorem}

\begin{proof}
Every finitely generated submodule $F$ has a composition series of length $r =
l_{A}(F)$. All its $r+1$ terms are submodules of $M$ of pairwise different
composition lengths. By the iso-multiplication property these yield $r+1$ distinct isomorphism
types among the ideal multiples $MB$; therefore $r \le u-1$. Choose a finitely
generated $F \le M$ of maximal possible length, which exists since the
possible lengths lie in the finite set $\{0, \dots, u-1\}$. If $F \neq M$,
select $m \in M \setminus F$. The finitely generated module $F + mA$ has
finite length, and the nonzero quotient $(F+mA)/F$ forces $l_{A}(F+mA) >
l_{A}(F)$, a contradiction. Hence $M = F$ and $l_{A}(M) \le u-1$.
Finite-length modules satisfy both classical chain conditions. The bound is
sharp, as witnessed by $A = k[x]/(x^{u-1})$ and $M = A_{A}$, which has $u$
distinct ideal multiples and length $u-1$ (for $u > 2$).
\end{proof}

\begin{remark}
The mere finiteness of
$\{[MB]\}$ does not imply finite length: for a principal right ideal domain
which is not a field, its regular module has only the two types $[0]$ and
$[A]$ among ideal multiples but may have infinite composition length. The
theorem uses finite length of finitely generated submodules, not just a finite
type spectrum.
\end{remark}

\begin{theorem}
\label{thm:artinian_bound} Let $A$ be a
right Artinian ring having exactly $t < \infty$ two-sided ideals, counted
including $0$ and $A$. If $M_{A}$ is an iso-multiplication module, then $M$ has finite
composition length, and
\[
l_{A}(M) \le t-1.
\]
In particular $M$ is finitely generated. The bound is best possible for each
integer $t \ge2$.
\end{theorem}

\begin{proof}
By the Hopkins-Levitzki theorem, $A$ is right Noetherian. Every finitely
generated right $A$-module consequently has finite composition length: indeed,
it is both Noetherian and Artinian. Let $F \le M$ be finitely generated, and
write $r = l_{A}(F)$. A composition series for $F$ is also a chain of
submodules of $M$. The iso-multiplication condition and Proposition \ref{prop:length_spectrum}%
's argument show that $r+1$ distinct isomorphism classes of ideal multiples
occur. Thus $r \le t-1$ for every finitely generated $F \le M$.

Among the finitely generated submodules choose one, say $F$, with maximal
possible length; its length exists because these lengths are integers in the
finite set $\{0, \dots, t-1\}$. If $F \neq M$, choose $x \in M \setminus F$.
Then $F+xA$ is finitely generated. The quotient $(F+xA)/F$ is nonzero and has
finite length, so $l_{A}(F+xA) > l_{A}(F)$, contrary to maximality. Therefore
$F = M$ and $l_{A}(M) \le t-1$.

For every $t \ge2$, take $A = k[x]/(x^{t-1})$ over a field $k$. Its exactly
$t$ two-sided ideals are the distinct ideals $(x^{j})$ for $0 \le j \le t-1$,
where $(x^{0}) = A$ and $(x^{t-1}) = 0$. Every submodule of $A_{A}$ is an
ideal, so $A_{A}$ is a classical multiplication module, hence iso-multiplication, and $l_{A}(A) = t-1$.
\end{proof}

\begin{corollary}
If $A$ is right
Artinian with finitely many two-sided ideals, every iso-multiplication module $M_{A}$ is a
finite-length module. In particular, it is Noetherian and Artinian and hence
satisfies both chain conditions up to isomorphism in the sense of Facchini-Nazemian.
\end{corollary}

\begin{proof}
Theorem \ref{thm:artinian_bound} yields finite composition length. The
remaining assertions are standard consequences of a composition series; no
cancellation or mutual-embedding hypothesis is needed.
\end{proof}

\begin{theorem}
\label{thm:equality_rigidity} Let $A$ be a
right Artinian ring with precisely $t < \infty$ two-sided ideals, and let
$M_{A}$ be an iso-multiplication module satisfying $l_{A}(M) = t-1$. Then the lengths of the
ideal multiples $MB$ are exactly $0, 1, \dots, t-1$, and exactly one two-sided
ideal $B$ gives each length. Consequently, if $N, L \le M$ then
\[
l_{A}(N) = l_{A}(L) \implies N \cong_{A} L.
\]
In particular, composition length induces a bijection from the set of
isomorphism classes of submodules of $M$ onto $\{0, \dots, t-1\}$. The
conclusion does not assert that $\mathcal{L}(M)$ is linearly ordered.
\end{theorem}

\begin{proof}
Put $r = l_{A}(M) = t-1$, and fix a composition series $0 = M_{0} \subsetneq
M_{1} \subsetneq\dots\subsetneq M_{r} = M$. For each $j$, the iso-multiplication property
gives a two-sided ideal $B_{j}$ satisfying $M_{j} \cong_{A} MB_{j}$, hence
$l_{A}(MB_{j}) = j$. Since these lengths are different, $B_{0}, \dots, B_{r}$
are $r+1 = t$ distinct two-sided ideals. They exhaust all the ideals of $A$,
and their corresponding ideal multiples have pairwise distinct lengths. Every
submodule $N \le M$ has some length $j \in\{0, \dots, r\}$ and is isomorphic
to an ideal multiple $MB$. Necessarily $B = B_{j}$, so $N \cong_{A} MB_{j}$.
The same reasoning applied to $L$ proves the assertion for equal lengths.
Every length $j$ occurs along the chosen composition series, proving the
bijection. No step implies comparability of different submodules.
\end{proof}

\begin{remark}
For $A = k[x]/(x^{r})$ and $M = A_{A}$, the $r+1$ ideals
are $(1), (x), \dots, (x^{r}) = 0$ and $l_{A}(M) = r$, so equality occurs for
every $r \ge1$. This familiar uniserial example verifies sharpness, but the
conclusion of Theorem \ref{thm:equality_rigidity} is uniqueness of isomorphism
type by length, not uniseriality in general. More abstractly, the same proof
works when $M$ has finite length $r$ and there are precisely $r+1$ distinct
isomorphism types among the ideal multiples $MB$: the $r+1$ distinct lengths
supplied by a composition series then exhaust those types.
\end{remark}

\begin{proposition}
\label{prop:recognition_test} Let $M_{A}$
have finite length $r$. Assume that any two submodules of $M$ having the same
composition length are isomorphic and that for every $j \in\{0, \dots, r\}$
some two-sided ideal $B_{j}$ satisfies $l_{A}(MB_{j}) = j$. Then $M$ is an
iso-multiplication module. Conversely, an iso-multiplication module of finite length always satisfies the
second assumption, whether or not the first assumption holds.
\end{proposition}

\begin{proof}
For $N \le M$ set $j = l_{A}(N)$. By assumption, $l_{A}(MB_{j}) = j$ and hence
$N \cong_{A} MB_{j}$. For the converse, take a composition series of $M$ and
use the iso-multiplication property on its terms. This establishes the length realization, but
does not automatically yield uniqueness of submodule types of equal length.
\end{proof}

\begin{example}
Let $k$ be a field and $A = k
\times k$. For $M = k \oplus k$ with the natural coordinatewise action the two
simple submodule types are different, but $M$ is a classical multiplication
module (all four submodules are ideal multiples). Thus the first hypothesis of
Proposition \ref{prop:recognition_test} is sufficient, not necessary:
iso-multiplication does not force equal-length submodules to be isomorphic in general.
\end{example}

\begin{example}
The algebra $T_{2}(k)$
has five two-sided ideals, so Theorem \ref{thm:artinian_bound} gives $l_{A}(M)
\le4$ for every iso-multiplication module $M$. The stronger complete classification in Theorem
\ref{thm:T2_classification} shows that the largest possible length is actually
three, attained by $A_{A} \cong P \oplus S_{2}$. The counting theorem is thus
a universal constraint rather than a substitute for a detailed classification.
\end{example}

\begin{remark}
For a finite-length module set
$\Lambda_{A}(M) = \{l_{A}(MB) : B \triangleleft A\}$. An iso-multiplication module has $\{0,
\dots, l_{A}(M)\} \subseteq\Lambda_{A}(M)$, but the converse is not valid:
matching lengths or composition-factor multiplicities does not produce module
isomorphisms. An ideal-length profile is a useful obstruction, not a complete characterization.
\end{remark}

\subsection{Composition-factor vectors}

Composition length does not distinguish the simple composition factors of
submodules. The following proposition gives a necessary condition in terms
of Jordan--H\"older multiplicities.

\begin{proposition}
\label{prop:factor_vector_obstruction}
Suppose $M_{A}$ has finite composition length, and let $S_{1}, \dots, S_{q}$
represent the isomorphism classes of simple composition factors of $M$. For $X
\le M$ put
\[
c(X) = ([X : S_{1}], \dots, [X : S_{q}]) \in\mathbb{N}^{q},
\]
where $[X : S_{i}]$ is the Jordan--H\"older multiplicity of $S_{i}$. If
$M$ is an iso-multiplication module, then
\[
\{c(N) : N \le M\} = \{c(MB) : B \triangleleft A\}.
\]
Thus a single composition-factor vector of a submodule which is not realized
by any ideal multiple disproves iso-multiplication.
\end{proposition}

\begin{proof}
For each $N \le M$, the iso-multiplication property yields $N \cong_{A} MB$ for some two-sided
ideal $B$. Isomorphic finite-length modules have the same Jordan--H\"older
multiplicities, so $c(N) = c(MB)$. Conversely, $MB \le M$ for every two-sided
ideal $B$, proving the reverse inclusion. The equality follows.
\end{proof}

\begin{example}
Let $A = T_{2}(k)$,
with radical $J = ke_{12}$, and two-sided ideals $L = ke_{11} \oplus ke_{12}$
and $K = ke_{12} \oplus ke_{22}$. Write $S_{1} = e_{11}A/e_{12}A$ and $S_{2} =
e_{22}A$ for the two nonisomorphic simple right modules, and set
\[
M = S_{1} \oplus S_{1} \oplus S_{2}.
\]
The five two-sided ideals of $A$ are $0, J, L, K, A$. Their actions on $M$
are
\[
M0 = MJ = 0, \quad ML = S_{1}^{\oplus2}, \quad MK = S_{2}, \quad MA = M.
\]
Thus the lengths of ideal multiples realize every value in $\{0,1,2,3\}$,
which is the entire possible scalar length spectrum of the submodules of $M$.
Nevertheless $S_{1} \le M$ has factor vector $(1,0)$, while the ideal
multiples have only the vectors
\[
(0,0), \quad(2,0), \quad(0,1), \quad(2,1).
\]
Therefore $S_{1}$ is not isomorphic to any $MB$, and $M$ is not iso-multiplication. In
particular, merely realizing all intermediate composition lengths by ideal
multiples is not sufficient for the defining property.
\end{example}

\begin{remark}
When a simple submodule $S \le M$ occurs, the
iso-multiplication property requires an ideal $B$ for which $MB \cong_{A} S$, not merely
$l_{A}(MB) = 1$. In a finite-length module the vector $c(S)$ has precisely one
nonzero coordinate, equal to 1. Therefore failure of any ideal multiple to
realize that coordinate vector excludes the iso-multiplication property. This is a practical
necessary test, not a sufficient criterion: distinct nonsplit extensions may
have the same Jordan--H\"older factor vector.
\end{remark}

\begin{remark}
The scalar-length and factor-vector tests are
obstructions, not general recognition theorems. Even matching complete
composition-factor vectors is weaker than matching module isomorphism classes
in nonsemisimple categories. Proposition \ref{prop:recognition_test} does give
a sufficient criterion, but only because it adds the independent hypothesis
that submodules with the same length are isomorphic.
\end{remark}

\subsection{General consequences}

\begin{theorem}
\label{thm:simple_rings} If $A$ is a simple ring and $M_{A}
\neq0$, then $M$ is iso-multiplication if and only if every nonzero submodule $N$ of $M$ is
isomorphic to $M$.
\end{theorem}

\begin{proof}
The only two-sided ideals are $0, A$, so the only ideal multiples are $0, M$.
\end{proof}

\begin{corollary}
For a simple Artinian ring $A$, the iso-multiplication modules are precisely the zero module
and the simple right modules.
\end{corollary}

\begin{proof}
A right module over $A$ is semisimple; if it has more than one simple summand,
a simple submodule is not isomorphic to the whole module, contradicting
Theorem \ref{thm:simple_rings}.
\end{proof}

\begin{theorem}
\label{thm:semisimple_rigidity}
Let $A$ be any associative ring with identity and let $M_{A}$ be a semisimple
right module, not necessarily finitely generated. Then
\[
M \text{ is im} \iff M \text{ is m}.
\]
In particular, every simple isomorphism type occurs with multiplicity at most
one in an iso-multiplication semisimple module.
\end{theorem}

\begin{proof}
Write $M = \bigoplus_{\lambda\in\Lambda} S_{\lambda}^{(\kappa_{\lambda})}$,
with the $S_{\lambda}$ pairwise nonisomorphic simple modules and
$\kappa_{\lambda}$ cardinals. For every two-sided ideal $B$, simplicity
implies $S_{\lambda}B = 0$ or $S_{\lambda}B = S_{\lambda}$; consequently $MB = \bigoplus_{S_{\lambda}B = S_{\lambda}} S_{\lambda}^{(\kappa_{\lambda})}.$
Assume $M$ is iso-multiplication and fix $\mu$ with $\kappa_{\mu}> 0$. A single copy $S_{\mu}$
is a submodule and so is isomorphic to some $MB$. Comparing simple summand
multiplicities in the displayed decomposition forces $\kappa_{\mu}= 1$. Thus
$M$ is multiplicity-free. For any submodule \(N\leq M\), semisimplicity implies that \(N\) is a direct sum of simple modules. Since \(M\) contains at most one copy of each simple isomorphism type, every simple submodule of \(M\) is contained in a unique simple direct summand \(S_\lambda\). Indeed, the coordinate projections of a simple submodule into the summands \(S_\lambda\) are homomorphisms, and a nonzero homomorphism between simple modules is an isomorphism. Since the \(S_\lambda\) are pairwise nonisomorphic, at most one coordinate projection can be nonzero. Consequently, every submodule \(N\) is the direct sum of a uniquely determined subset of the simple direct summands of \(M\). If \(N\cong M B\), comparison of the multiplicities of the simple isomorphism types shows that this subset is exactly the set of summands occurring in \(M B\). Hence \(N=M B\), and \(M\) is an \(m\)-module. The converse is immediate.
\end{proof}

\begin{corollary}
If $A$ is semisimple Artinian, then $M_{A}$ is iso-multiplication
if and only if it is a multiplication module, precisely when each simple $A$-module occurs with
multiplicity at most one in $M$. In particular, this applies without assuming
$M$ finitely generated.
\end{corollary}

\begin{corollary}
If a semisimple module $M$ contains two
isomorphic nonzero simple direct summands, then $M$ is not iso-multiplication. This applies,
in particular, to $S \oplus S$ for every simple $S$ over every ring.
\end{corollary}

\begin{proposition}
For every right $A$-module $M$, set
\[
\mathcal{I}(M) = \{[MB] : B \triangleleft A\}, \quad\mathcal{S}(M) = \{[N] : N
\in\mathcal{L}(M)\},
\]
where brackets denote right-module isomorphism types. Then $M$ is iso-multiplication if and
only if $\mathcal{I}(M) = \mathcal{S}(M)$. Consequently, if $A$ has only $t$
two-sided ideals and $M$ is iso-multiplication, at most $t$ isomorphism types occur among the
submodules of $M$.
\end{proposition}

\begin{proof}
Every $MB$ is a submodule, so $\mathcal{I}(M) \subseteq\mathcal{S}(M)$ without
assumptions. The reverse inclusion is exactly the definition of iso-multiplication. The
cardinality estimate is immediate.
\end{proof}

\begin{remark}
The preceding equality gives an obstruction to the iso-multiplication
property. It does not characterize rings admitting iso-multiplication modules.
\end{remark}

\begin{theorem}
\label{thm:finite_direct_products} Let $A = A_{1}
\times\dots\times A_{n}$ and $M = \bigoplus_{i=1}^{n} M_{i}$, where $M_{i} =
Me_{i}$ is regarded as a right $A_{i}$-module. Then
\[
M_{A} \text{ is im} \iff(M_{i})_{A_{i}} \text{ is im for every } i.
\]

\end{theorem}

\begin{proof}
Every right submodule $N$ decomposes as $\bigoplus N e_{i}$, and each
two-sided ideal as $B_{1} \times\dots\times B_{n}$. Consequently $MB =
\bigoplus M_{i} B_{i}$. Componentwise isomorphisms prove sufficiency.
Conversely, given $N_{i} \le M_{i}$, view it as a submodule of $M$ supported
only in component $i$. An isomorphism $N_{i} \cong MB$ must preserve each
central idempotent and hence restrict to $N_{i} \cong M_{i} B_{i}$, proving necessity.
\end{proof}

\begin{proposition}
Suppose $M_{A}$ has finite
composition length. If $M$ is iso-multiplication, then the multiplicity vector of the simple
composition factors of each submodule $N$ coincides with that of $MB$ for some
two-sided ideal $B$. Equality of lengths or composition-factor vectors alone
is not a sufficient criterion for the iso-multiplication property.
\end{proposition}

\begin{proof}
The first assertion follows from an isomorphism $N \cong MB$ and the
Jordan--H\"older theorem. For the warning, nonisomorphic finite-length
modules can have the same composition factors (for example $P$ and $S_{1}
\oplus S_{2}$ over $T_{2}(k)$).
\end{proof}

\subsection{Local Loewy rigidity}

We now consider local rings. The next theorem holds without a commutativity
assumption; its converse fails for noncommutative rings.

\begin{theorem}
\label{thm:local_loewy} Let $A$ be any
(not necessarily commutative) local ring with Jacobson radical $J$, and let
$M_{A}$ be a right module satisfying $MJ^{h} = 0$ for some integer $h \ge1$.
If $M$ is an iso-multiplication module, then $M$ is cyclic. Neither $A$ nor
$M$ is assumed to be noetherian, artinian, or finitely generated. In
particular, every iso-multiplication module over a right artinian local ring
is cyclic and has finite composition length.
\end{theorem}

\begin{proof}
Assume $M \neq0$ and choose $h$ minimal with $MJ^{h} = 0$. Then $MJ^{h-1}
\neq0$, so there is an $x \in M$ with $xJ^{h-1} \neq0$. As $xA \le M$ and $M$
is iso-multiplication, $xA \cong_{A} MB$ for some two-sided ideal $B
\triangleleft A$. If $B \neq A$ then $B \subseteq J$ because $A$ is local:
every proper two-sided ideal is contained in its unique maximal right ideal.
Consequently
\[
(MB)J^{h-1} = M(BJ^{h-1}) \subseteq MJ^{h} = 0.
\]
On the other hand, $(xA)J^{h-1}$ contains $xJ^{h-1} \neq0$, contradicting the
right $A$-module isomorphism $xA \cong MB$. Thus $B = A$ and $M \cong xA$. An
image of a cyclic module under an isomorphism is cyclic; hence $M$ is cyclic.
If $A$ is right artinian, $J$ is nilpotent and a cyclic right $A$-module has
finite length.
\end{proof}

\begin{remark}
For commutative local $A$,
cyclic modules are ordinary multiplication modules, so Theorem
\ref{thm:local_loewy} recovers the full equivalence in the restricted
commutative comparison. For noncommutative local rings, only the implication
``iso-multiplication $\implies$ cyclic'' has been proved. The reverse
implication is false, as the following example demonstrates. In particular, it
would be an error to replace ``commutative'' by ``noncommutative'' in the
cyclicity characterizations of the commutative cases.
\end{remark}

\begin{example}
\label{ex:cyclic_not_im}
Let $k$ be any field, put $D = k(t)$, and let
\[
T =
\begin{pmatrix}
0 & t\\
1 & 0
\end{pmatrix}
\in M_{2}(D).
\]
The polynomial $X^{2}-t$ is irreducible over $D=k(t)$. Indeed, if
$t=(f(t)/g(t))^{2}$ with nonzero coprime polynomials $f,g\in k[t]$, then
$f(t)^{2}=t g(t)^{2}$. Since $t$ is prime in $k[t]$, it divides $f$, say
$f=th$. Substitution gives $t^{2}h^{2}=t g^{2}$, hence $g^{2}=t h^{2}$ and $t\mid g^{2}$, so
$t\mid g$. This contradicts $\gcd(f,g)=1$. Thus $t$ is not a square in $D$ and
the quadratic polynomial $X^{2}-t$ is irreducible. Thus $T$ has no
one-dimensional invariant $D$-subspace in $D^{2}$. The matrix $T$ is
transcendental over $k$: if $p(T) = 0$ for $p \in k[X]$, write $p(X) =
p_{0}(X^{2}) + X p_{1}(X^{2})$. As $T^{2} = t I_{2}$, we have $p(T) =
p_{0}(t)I_{2} + p_{1}(t)T$; the off-diagonal entries imply $p_{1}(t) = 0$,
then the diagonal entries imply $p_{0}(t) = 0$. As $t$ is transcendental over
$k$, both polynomials vanish, so $p = 0$. Thus the rule $t \mapsto T$ extends
to a unital field embedding $\rho: D = k(t) \hookrightarrow M_{2}(D)$,
because, writing $p(T) = u(t)I_{2} + v(t)T$, its determinant is $u(t)^{2} - t
v(t)^{2}$, which is nonzero whenever $(u,v) \neq(0,0)$ by the irreducibility
of $X^{2} - t$ over $D$. Take $V = D^{2}$ as a right $D$-vector space of
column vectors, with coordinatewise right multiplication $v \cdot b = (v_{1}
b, v_{2} b)^{T}$. Independently define its left action by $a \cdot v =
\rho(a)v$. The matrices $\rho(a)$ act $D$-linearly on the right, so $a
\cdot(vb) = (a \cdot v)b$ as required for a bimodule. The two actions do not commute in the sense required for a commutative
trivial extension: for $v=(1,0)^T$, one has
$\rho(t)v=(0,1)^T\neq(t,0)^T=vt$. Thus the ring constructed below is
noncommutative. This makes $V$ a $(D,D)$-bimodule with no nonzero proper
subbimodule: a nonzero proper subbimodule would have to be a right $D$-line
invariant under $\rho(t) = T$.
Such a line would supply an eigenvector over $D$, contradicting the
irreducibility of $X^{2} - t$.

Form the trivial extension
\[
A = D \times V, \quad(a,v)(b,w) = (ab, a \cdot w + vb).
\]
The subset $J = 0 \times V$ satisfies $J^{2} = 0$, and $A/J \cong D$. An
element $(a,v)$ with $a \neq0$ is invertible: since $(0,V)^{2} = 0$, the usual
inverse of the unit $(a,0)$ extends to $(a,v)$ by a one-step nilpotent
correction. Conversely, elements of $J$ are nonunits. Thus $A$ is a local
noncommutative ring with Jacobson radical $J$. Its two-sided ideals are
exactly $0, J, A$: every proper two-sided ideal is contained in $J$, and such an ideal
is exactly a $(D,D)$-subbimodule of $V$ under the identification
$J=0\times V$.

As a right $A$-module, $J$ is a two-dimensional right $D$-vector space,
whereas every one-dimensional right $D$-subspace $W \le V$ gives a nonzero
right ideal $0 \times W$ of $A$. The ideal multiples of the regular right
module $A_{A}$ are exactly $0, J, A$. Their composition lengths are $0, 2, 3$,
respectively: $J$ is a direct sum of two copies of the unique simple
right $A$-module $D_A$, and $A/J\cong D_A$. The right ideal $0\times W$
has composition length 1.
Therefore $A_{A}$ is cyclic, but not iso-multiplication. This shows that the
converse to Theorem \ref{thm:local_loewy} fails even for a local ring with
square-zero Jacobson radical.
\end{example}

\begin{corollary}
For every right artinian local ring $A$,
every iso-multiplication right module is a finitely generated cyclic module.
If $A$ is commutative, the converse holds; for noncommutative $A$ it can fail
even for the regular module $A_{A}$.
\end{corollary}

\begin{proof}
The forward assertion follows from Theorem \ref{thm:local_loewy}. For
commutative $A$, the converse follows because every cyclic module is an
ordinary multiplication module. The explicit local ring in Example
\ref{ex:cyclic_not_im} has square-zero radical and is right artinian (it has a
finite composition series of length three as a right module over itself),
providing the claimed counterexample. Notice that ``right artinian'' here is a
property of $A_{A}$ and does not require finite dimension over the central
field $k$.
\end{proof}

The preceding examples separate multiplication modules from
iso-multiplication modules. The next section gives restrictions of a different
kind using Goldie localization.

\begin{remark}
For a division ring $D$ and $n > 2$, the
regular right module of $M_{n}(D)$ is not iso-multiplication. Its only two-sided ideals are
$0$ and $M_{n}(D)$, whereas the nonzero proper right ideal $e_{11}M_{n}(D)$
has length one and the regular module has length $n$. This elementary
observation foreshadows the general Goldie obstruction below.
\end{remark}

\section{Localization over semiprime right Goldie rings}

\label{sec:rational_rigidity}

We use the classical quotient ring of a semiprime right Goldie ring. The
necessary facts on Ore localization and Goldie's theorem can be found in
\cite{SmallNotes}. No commutativity assumption on products of right ideals
is made.

Let $A$ be semiprime right Goldie and $C$ the set of regular elements of $A$.
Then $C$ is a right denominator set and $Q = AC^{-1}$ is semisimple Artinian.
Write
\[
Q \cong\prod_{i=1}^{s} M_{n_{i}}(D_{i}),
\]
where the $D_{i}$ are division rings and $n_{i} \ge1$. For a right $A$-module
$M$, write $M_{Q} = M \otimes_{A} Q$. The module $M_{Q}$ is semisimple, with
possibly infinite multiplicities.

\begin{remark}
The $-\otimes_{A} Q$ is exact for the classical
Ore quotient used here, and two-sided ideals of the semisimple Artinian ring
$Q$ are sums of full simple factors. Neither statement is new. The substantive
step below is to apply these facts to the ideal multiples representing every
submodule, forcing multiplicity one. The resulting theorem is a necessary
condition and must not be interpreted as a recognition criterion.
\end{remark}

\begin{lemma}
\label{lem:localization_ideals} If $B$ is a
two-sided ideal of $A$, then $BQ$ is a two-sided ideal of $Q$. Moreover, for
every right $A$-module $M$ there is a natural $Q$-module isomorphism
\[
(MB) \otimes_{A} Q \cong M_{Q} (BQ).
\]
Consequently, $M_{Q} (BQ)$ is the direct sum of entire isotypic components of
$M_{Q}$.
\end{lemma}

\begin{proof}
We use the right Ore convention $Q = AC^{-1}$: elements are represented by
right fractions $ac^{-1}$ with $a \in A$ and $c \in C$. Put $BQ = \{ \sum_{j}
b_{j} q_{j} : b_{j} \in B, q_{j} \in Q \}$, a right $Q$-ideal. Because $B$ is
a two-sided ideal, $a(BQ) \subseteq BQ$ for every $a \in A$. For a regular
element $c \in C$, left multiplication
\[
L_{c} : BQ \longrightarrow BQ, \quad y \mapsto cy
\]
is right $Q$-linear and injective: $c$ is a unit in $Q$. As a right submodule
of the finite-length regular module $Q_{Q}$, the module $BQ$ has finite
length. Therefore the injective endomorphism $L_{c}$ is surjective. In
particular, for each $y \in BQ$ the unique solution to $cz = y$ in $Q$ belongs
to $BQ$, and hence $c^{-1}(BQ) \subseteq BQ$. Thus $(ac^{-1})(BQ) \subseteq
BQ$ for all $a \in A$ and $c \in C$: the right ideal $BQ$ is also a left
$Q$-ideal. Since $Q$ is semisimple Artinian, a central idempotent $e$
satisfies $BQ = eQ$.

For a classical right Ore localization at regular elements, $Q$ is flat as a
left $A$-module (this is the standard Ore-localization exactness theorem);
thus the $-\otimes_{A} Q$ is exact. More explicitly, apply it to $0 \to MB \to
M$; exactness gives an embedding $(MB) \otimes_{A} Q \hookrightarrow M_{Q}$.
The image is the right $Q$-span of elements $(mb) \otimes1$ with $m \in M$ and
$b \in B$. As $(mb) \otimes1 = (m \otimes1)b$, this image is exactly
$M_{Q}(BQ)$. This identifies the two localized modules without presupposing
that $MB \hookrightarrow M$ remains injective under an arbitrary nonflat
tensor product. Finally $M_{Q} e$ is the full direct sum of the simple
constituents supported on the simple factors of $Q$ on which $e$ is 1; it
cannot select only some copies of one simple type.
\end{proof}

\begin{lemma}
\label{lem:simple_constituent_detected} For every simple $Q$-submodule $S$ of
$M_{Q}$ there exists a cyclic right $A$-submodule $N$ of $M$ such that
\[
N \otimes_{A} Q \cong_{Q} S.
\]
This holds even if $M$ is not finitely generated or has nonzero $C$-torsion.
\end{lemma}

\begin{proof}
Choose $0 \neq x \in S$. Because $Q$ is a right Ore localization, write $x =
\sum_{j=1}^{t} (m_{j} \otimes1) q_{j}$ and put the finitely many $q_{j}$ over
a common right denominator $c \in C$: $x = (m \otimes1) c^{-1}$ for some $m
\in M$. Thus $m \otimes1 = xc$ is nonzero and belongs to $S$. For $N = mA$,
exact localization identifies $N \otimes_{A} Q$ with $(m \otimes1)Q$, which
equals $S$ by simplicity.
\end{proof}

\begin{theorem}
\label{thm:rational_multiplicity_free}
Let $A$ be an arbitrary semiprime right Goldie ring and $M_{A}$ an
iso-multiplication module. Then
\[
M \otimes_{A} Q_{\text{cl}}^{r}(A)
\]
is multiplicity-free semisimple. Equivalently, if $Q \cong\prod_{i=1}^{s}
M_{n_{i}}(D_{i})$ and $V_{i}$ denotes the unique simple right module for its
$i$th matrix factor, then
\[
M \otimes_{A} Q \cong_{Q} \bigoplus_{i \in E} V_{i}
\]
for some $E \subseteq\{1, \dots, s\}$. No finiteness hypothesis on $M$ is required.
\end{theorem}

\begin{proof}
Fix an index $i$ for which the multiplicity of $V_{i}$ in $M_{Q}$ is nonzero,
and choose a simple submodule $S \cong V_{i}$ of $M_{Q}$. By Lemma
\ref{lem:simple_constituent_detected}, some cyclic submodule $N \le M$ has
$N_{Q} \cong S$. Since $M$ is iso-multiplication, $N \cong_{A} MB$ for a
two-sided ideal $B \triangleleft A$. Localizing this isomorphism and applying
Lemma \ref{lem:localization_ideals}, we obtain
\[
V_{i} \cong_{Q} N_{Q} \cong_{Q} (MB)_{Q} \cong_{Q} M_{Q} (BQ).
\]
The module on the right is a direct sum of entire isotypic components of
$M_{Q}$. Since it is a single simple module, the multiplicity of $V_{i}$ in
$M_{Q}$ is exactly one. This holds for every occurring simple type.
\end{proof}

\begin{corollary}
\label{cor:division_ring_factors} If $A$ is semiprime right Goldie and $A_{A}$
is iso-multiplication, then
\[
Q_{\text{cl}}^{r}(A) \cong D_{1} \times\dots\times D_{s}
\]
for division rings $D_{1}, \dots, D_{s}$. In particular, the regular module of
an order in a semisimple Artinian ring with a genuine matrix factor $M_{n}%
(D)$, $n \ge2$, cannot be iso-multiplication.
\end{corollary}

\begin{proof}
For $M = A_{A}$ we have $M_{Q} = Q_{Q}$. As a right module over $M_{n}(D)$,
the regular module is the direct sum of $n$ copies of its unique simple
module. Theorem \ref{thm:rational_multiplicity_free} forces every $n_{i} = 1$.
\end{proof}

\begin{corollary}
If $A$ is semiprime right Goldie and the regular
module $A_{A}$ is iso-multiplication, then every idempotent of $A$ is central.
\end{corollary}

\begin{proof}
By Corollary \ref{cor:division_ring_factors}, the canonical embedding $A
\hookrightarrow Q_{\text{cl}}^{r}(A)$ has codomain a finite direct product of
division rings. Every idempotent in that product is a central coordinate
idempotent, since a division ring has only the idempotents $0$ and $1$. Thus,
for $e^{2} = e \in A$ and $a \in A$, the images of $ea$ and $ae$ in $Q$ agree;
injectivity gives $ea = ae$ in $A$.
\end{proof}

\begin{remark}
Centrality of idempotents is also a standard consequence of reducedness. The
corollary is recorded only to make the structural restriction explicit, not as
a separate novelty claim.
\end{remark}

\begin{corollary}
If $A$ is prime right Goldie and $M_{A}$
is iso-multiplication, then $M_{Q}$ is either zero or a simple right module
over the simple Artinian ring $Q$. In particular, if $A_{A}$ is
iso-multiplication, then $A$ is a domain and $Q$ is a division ring.
\end{corollary}

\begin{proof}
Goldie's theorem makes $Q$ simple Artinian; Theorem
\ref{thm:rational_multiplicity_free} gives the first statement. For the
regular module $Q_{Q}$, multiplicity-freeness forces $Q$ to be a division
ring. Since $A \hookrightarrow Q$, nonzero elements of $A$ have no zero product.
\end{proof}

\begin{theorem}
\label{thm:exact_uniform_dimension} Let $A$
be semiprime right Goldie, let $Q = Q_{\text{cl}}^{r}(A)$, and suppose that
the canonical map $M \to M_{Q} = M \otimes_{A} Q$ is injective. If $M_{Q}$ has
finite composition length over $Q$, then
\[
\text{u.dim}_{A}(M) = l_{Q}(M_{Q}).
\]
In particular, if $M$ is also an iso-multiplication module and $Q \cong\prod_{i=1}^{s}
M_{n_{i}}(D_{i})$, then
\[
\text{u.dim}_{A}(M) = l_{Q}(M_{Q}) = |\{i : \text{the } i\text{th simple type
occurs in } M_{Q}\}| \le s.
\]
Neither finite generation of $M$ nor noetherianity of $A$ is assumed.
\end{theorem}

\begin{proof}
Let $d = l_{Q}(M_{Q})$. If $N_{1} \oplus\dots\oplus N_{t} \le M$ is a direct
sum of nonzero submodules, every $(N_{j})_{Q}$ is nonzero because $M \to
M_{Q}$ is injective. Exactness of Ore localization gives a direct sum of $t$
nonzero $Q$-submodules of $M_{Q}$, hence $t \le d$. This proves $\text{u.dim}%
_{A}(M) \le d$.

Conversely, choose a decomposition $M_{Q} = S_{1} \oplus\dots\oplus S_{d}$
into simple right $Q$-modules. For each $j$ select $0 \neq x_{j} \in S_{j}$.
The common-denominator argument in Lemma \ref{lem:simple_constituent_detected}
writes $x_{j} = (m_{j} \otimes1)c_{j}^{-1}$ with $m_{j} \in M$, $c_{j} \in C$.
Since $m_{j} \otimes1 = x_{j} c_{j} \in S_{j}$ is nonzero, the submodule
$m_{j} A$ of $M$ is nonzero and its image in $M_{Q}$ lies in $S_{j}$. The sum
$m_{1} A \oplus\dots\oplus m_{d} A$ is therefore direct inside $M$ (by
injectivity of $M \to M_{Q}$). Thus $\text{u.dim}_{A}(M) \ge d$. Equality
follows. The last displayed formula uses Theorem
\ref{thm:rational_multiplicity_free}.
\end{proof}

\begin{corollary}
\label{cor:faithful_goldie_saturation} Let $A$ be a semiprime right Goldie
ring, let $C$ be its set of regular elements, and write
\[
Q = Q_{\text{cl}}^{r}(A) \cong\prod_{i=1}^{s} M_{n_{i}}(D_{i}).
\]
Suppose that the right $A$-module $M$ is faithful, $C$-torsion-free, and
iso-multiplication. If $V_{i}$ denotes the unique simple right module
associated with the $i$th simple Artinian factor, then
\[
M \otimes_{A} Q \cong_{Q} \bigoplus_{i=1}^{s} V_{i}, \quad\text{u.dim}_{A}(M)
= s.
\]
In particular, no simple quotient-ring component is missing and no simple type
occurs twice. No finite-generation hypothesis on $M$ is needed.
\end{corollary}

\begin{proof}
Write $M_{Q} = M \otimes_{A} Q$. Since $M$ is $C$-torsion-free, the canonical
map $M \to M_{Q}$ is injective. By Theorem
\ref{thm:rational_multiplicity_free}, each simple right $Q$-module type occurs
in $M_{Q}$ with multiplicity at most one. We prove that every type occurs.

Fix $1 \le i \le s$ and let $e_{i}$ be the nonzero central idempotent of $Q$
corresponding to its $i$th factor. In the classical right Ore quotient ring,
express $e_{i} = ac^{-1}$ with $a \in A$ and $c \in C$. Thus $a = e_{i} c
\neq0$ in $Q$, hence also in $A$ since $A \hookrightarrow Q$ is injective;
moreover $a \in e_{i} Q \cap A$. This explicitly exhibits a nonzero element of
$A$ supported on the $i$th central component of $Q$. If $M_{Q} e_{i} = 0$,
then $(m \otimes1)a = (m \otimes1)e_{i} c = 0$ for every $m \in M$.
Injectivity of $M \to M_{Q}$ implies $ma = 0$ for all $m$, contradicting
faithfulness and $a \neq0$. Therefore $M_{Q} e_{i} \neq0$, so its $i$th simple
type occurs. Multiplicity-freeness now yields the asserted direct sum. Finally
Theorem \ref{thm:exact_uniform_dimension} gives $\text{u.dim}_{A}(M) =
l_{Q}(M_{Q}) = s$.
\end{proof}

\begin{remark}
This is a consequence of the
rational rigidity and uniform-dimension theorems, not an independent
necessary-and-sufficient characterization. The roles of its two extra
assumptions differ: faithfulness prevents a central quotient-ring component
from disappearing, while $C$-torsion-freeness permits detection of the
annihilation relation in the original module.
\end{remark}

\begin{example}
For $A = \mathbb{Z} \times\mathbb{Z}$ and $M =
A_{A}$, the iso-multiplication property is classical multiplication and $Q = \mathbb{Q}
\times\mathbb{Q}$. We have $\text{u.dim}_{A}(M) = 2 = l_{Q}(M_{Q})$, so the
general upper bound $s$ can be attained. For the skew polynomial right
principal ideal domain of Example \ref{ex:skew_polynomial}, the corresponding
values are both 1.
\end{example}

\begin{corollary}
Let $A$ be prime right Goldie, and let $M
\neq0$ be a $C$-torsion-free iso-multiplication right $A$-module. Then $M$ is
uniform. In particular, every nonzero pair of submodules of $M$ has nonzero intersection.
\end{corollary}

\begin{proof}
Here $Q$ is simple Artinian. Theorem \ref{thm:rational_multiplicity_free} says
$M_{Q}$ is either zero or simple; torsion-freeness and $M \neq0$ rule out
zero. If $N_{1}, N_{2} \le M$ are nonzero with zero intersection, then $N_{1}
\oplus N_{2}$ embeds in $M$. Localizing gives two nonzero $Q$-submodules in
direct sum inside the simple module $M_{Q}$, a contradiction.
\end{proof}

\begin{corollary}
Let $M_{A}$ be iso-multiplication and put
$\overline{A} = A / \text{Ann}_{A}(M)$. If $\overline{A}$ is semiprime right
Goldie, rational rigidity applies to $M$ as a right $\overline{A}$-module. If
$M$ is also torsion-free over the regular elements of $\overline{A}$, the
exact uniform-dimension formula applies.
\end{corollary}

\begin{proof}
Every submodule of $M$ is the same over $A$ and $\overline{A}$. If $N \cong
MB$ as $A$-modules for a two-sided ideal $B \triangleleft A$, this is an
$\overline{A}$-isomorphism and $MB = M\overline{B}$ for $\overline{B} = (B +
\text{Ann}_{A}(M)) / \text{Ann}_{A}(M) \triangleleft\overline{A}$. Hence $M$
is iso-multiplication over $\overline{A}$, and the assertions follow from
Theorems \ref{thm:rational_multiplicity_free} and
\ref{thm:exact_uniform_dimension}.
\end{proof}

\begin{remark}
The uniform-dimension conclusions do
not follow merely from $M_{Q} = 0$: a torsion module may disappear under
localization. Theorem \ref{thm:exact_uniform_dimension} deliberately requires
the canonical map $M \to M_{Q}$ to be injective; this is not automatic for
iso-multiplication modules.
\end{remark}

\begin{example}
The skew polynomial ring
$A = k(t)[x;\sigma]$ of Example \ref{ex:skew_polynomial} is a noncommutative
principal right ideal domain. Its right classical quotient ring is a division
ring, and $A_{A}$ is iso-multiplication without being a classical multiplication module.
Thus Corollary \ref{cor:division_ring_factors} does not force commutativity,
right duo behavior, or the ordinary multiplication property.
\end{example}

\begin{example}
The right Artinian algebra
$A = T_{2}(k)$ is right Goldie, and its regular module is iso-multiplication
by Example \ref{ex:basic_counterexample}. Yet $A$ is not semiprime: its
radical $ke_{12}$ has square zero. Thus the semiprime hypothesis is essential
for applying the above quotient-ring mechanism; the example does not provide a
semiprime counterexample to the conclusion.
\end{example}

\begin{remark}
A zero rational localization is vacuously
multiplicity-free and therefore cannot characterize the iso-multiplication property; for
example $(\mathbb{Z}/p\mathbb{Z})^{2}$ is not iso-multiplication over $\mathbb{Z}$ although
its tensor product with $\mathbb{Q}$ vanishes. The next example shows that
even a nonzero simple rational localization is insufficient.
\end{remark}

\begin{example}
\label{ex:rational_not_sufficient} Let $A = k[x]$, $S = A/(x)$, and $M = A
\oplus S$. With $Q = k(x)$ we have
\[
M_{Q} \cong Q,
\]
so the localization is nonzero and simple. Nevertheless $M$ is not an
iso-multiplication module. Indeed, for each ideal $B = (f)$ the isomorphism type of $MB$ is
$0$ if $f=0$, is $[A]$ if $0 \neq f \in(x)$, and is $[A \oplus S]$ if $f
\notin(x)$. But $S \le M$ is not isomorphic to any of these. Hence rational
multiplicity-freeness is strictly a necessary condition; even a simple nonzero
localization is not enough. This contrasts with the preceding example where
all of $M$ vanished under localization.
\end{example}

\begin{remark}
Corollary
\ref{cor:division_ring_factors} strengthens Lemma
\ref{lem:matrix_size_obstruction} in a different direction: it applies to
every semiprime right Goldie ring, without any principal-ideal,
right-noetherian, or isoartinian assumption. It is a deduction from the stated
Goldie rigidity theorem, not a new assertion of Goldie's theorem. The conclusion is obtained from the Goldie rigidity theorem.
\end{remark}

\subsection{Rational submodules and essentiality}

We also describe the submodules of $M_Q$ and the essential submodules
detected by localization. We retain the notation $C,Q$ from Section
\ref{sec:rational_rigidity} and denote the image of $N\le M$ in $M_Q$ by
$N_Q$.

\begin{lemma}
\label{lem:rational_submodule_lift} Let $A$ be semiprime right Goldie and
$M_{A}$ any module. For each $Q$-submodule $W \le M_{Q}$ there is a submodule
$N \le M$ satisfying $N_{Q} = W$. If $M \to M_{Q}$ is injective, one may take
\[
N = \{m \in M : m \otimes1 \in W\}.
\]

\end{lemma}

\begin{proof}
Let $N$ be the preimage of $W$ under the canonical map $M \to M_{Q}$. Clearly
$N_{Q} \subseteq W$. Conversely, every $w \in W$ can be written as $(m
\otimes1)c^{-1}$ with $m \in M$ and $c \in C$ by clearing finitely many right
denominators. Then $m \otimes1 = wc \in W$, so $m \in N$ and
$w=(m\otimes1)c^{-1}\in N_Q$. Thus $N_Q=W$. The same inverse-image
formula defines $N$ even if $M\to M_Q$ is not injective; under injectivity
it may additionally be identified with the intersection of $M$ and $W$
inside $M_Q$.
\end{proof}

\begin{theorem}
\label{thm:boolean_localization} Let $A$ be
semiprime right Goldie and let $M_{A}$ be iso-multiplication. Write
\[
M_{Q} \cong\bigoplus_{i \in E} V_{i}
\]
with pairwise nonisomorphic simple right $Q$-modules $V_{i}$, as guaranteed by
Theorem \ref{thm:rational_multiplicity_free}. Then the map
\[
\mathcal{L}(M) \longrightarrow\mathcal{L}_{Q}(M_{Q}), \quad N \mapsto N_{Q}
\]
is surjective, and $\mathcal{L}_{Q}(M_{Q})$ is a finite Boolean lattice with
$2^{|E|}$ elements. Moreover every localized submodule is a localized ideal
multiple: for each $N \le M$ there is $B \triangleleft A$ such that
\[
N_{Q} = (MB)_{Q} = M_{Q} (BQ).
\]
Thus each rational submodule is of the form $M_{Q} e$ for a central idempotent
$e \in Q$; $e$ is unique up to components annihilating $M_{Q}$. In particular,
the isomorphism-class information of localized submodules is controlled by the
support set $E$.
\end{theorem}

\begin{proof}
Surjectivity follows from Lemma \ref{lem:rational_submodule_lift}. Distinct
$V_{i}$ are nonisomorphic simple modules, so each $Q$-submodule of their
direct sum is uniquely the sum of a subset of the $V_{i}$. This gives the
Boolean lattice. By the iso-multiplication property, $N \cong_{A} MB$ for some two-sided ideal
$B$; hence $N_{Q} \cong_{Q} (MB)_{Q}$. Lemma \ref{lem:localization_ideals}
shows that $(MB)_{Q}$ is a sum of full isotypic components. Since $M_{Q}$ is
multiplicity-free, any two isomorphic $Q$-submodules of $M_{Q}$ are equal,
proving the asserted equality inside $M_{Q}$.
\end{proof}

\begin{theorem}
\label{thm:essentiality_rational}
Let $A$ be semiprime right Goldie and suppose $M \to M_{Q}$ is injective. For
every $N \le M$,
\[
N \text{ is essential in } M \iff N_{Q} = M_{Q}.
\]
If, in addition, $M$ is iso-multiplication, this is equivalent to saying that
$N_{Q}$ contains every simple $Q$-type occurring in $M_{Q}$.
\end{theorem}

\begin{proof}
First suppose $N_{Q}=M_{Q}$. Let $0\neq L\leq M$. The injectivity of
$M\longrightarrow M_{Q}$ implies $L_{Q}\neq0$. Because Ore localization
is exact, it preserves the intersection of the embedded submodules $N,L$:
\[
(N\cap L)_{Q}=N_{Q}\cap L_{Q}=M_{Q}\cap L_{Q}=L_{Q}\neq0.
\]
Thus $N\cap L\neq0$, proving that $N$ is essential in $M$.

For the reverse implication, suppose that $N$ is essential in $M$ but
$N_{Q}\neq M_{Q}$. Since $Q$ is semisimple Artinian, $N_Q$ is a direct
summand of $M_Q$. Choose a nonzero simple submodule $S$ of a complementary
summand; then $S\cap N_Q=0$. Lemma
\ref{lem:simple_constituent_detected} supplies a cyclic submodule
$L=mA\leq M$ whose localized image is $L_Q=S$. The equality $L_Q=S\neq0$
shows that $L\neq0$. Because the localization functor is exact, it
preserves finite intersections of submodules inside $M$; explicitly,
\[
(N\cap L)_{Q}=N_{Q}\cap L_{Q}=N_{Q}\cap S=0.
\]
As $M\to M_Q$ is injective, its restriction to $N\cap L$ is also
injective. Thus a nonzero element of $N\cap L$ would have a nonzero image
in $(N\cap L)_Q$, contradicting the displayed equality. Hence
$N\cap L=0$, in conflict with the essentiality of $N$. This proves
$N_Q=M_Q$.
If $M$ is also iso-multiplication, Theorem
\ref{thm:rational_multiplicity_free} decomposes $M_{Q}$ into pairwise
nonisomorphic simple constituents; thus $N_{Q}=M_{Q}$ is equivalent to
$N_{Q}$ containing all their types.
\end{proof}

\begin{corollary}
Under the hypotheses of Theorem
\ref{thm:boolean_localization}, the ring $\text{End}_{Q}(M_{Q})$ is a finite
direct product of division rings (with the zero ring in case $M_{Q} = 0$).
\end{corollary}

\begin{proof}
The $V_{i}$ are pairwise nonisomorphic simple modules. Therefore all
off-diagonal Hom groups vanish and
\[
\text{End}_{Q}(M_{Q}) \cong\prod_{i \in E} \text{End}_{Q}(V_{i}),
\]
where each factor is a division ring by Schur's lemma. Since $E$ is
a subset of the finite set of simple factors of $Q$, the product is finite.
\end{proof}

\begin{example}
Let $A = \mathbb{Z} \times\mathbb{Z}$ and
$M = A_{A}$. Then $M_{Q} = \mathbb{Q} \oplus\mathbb{Q}$ with two nonisomorphic
simple modules over $Q = \mathbb{Q} \times\mathbb{Q}$. Its localized submodule
lattice is the four-element Boolean lattice. The submodule $N = (2\mathbb{Z})
\times(3\mathbb{Z})$ is essential, and $N_{Q} = M_{Q}$; by contrast,
$\mathbb{Z} \times0$ is not essential, and its localization is only one
Boolean atom.
\end{example}

\begin{remark}
The Boolean conclusion is a statement about the
localized lattice, not a claim that $\mathcal{L}(M)$ itself is Boolean. For
example, the regular module over $\mathbb{Z}$ is iso-multiplication but its
actual submodule lattice contains many proper nonzero ideals. The torsion-free
hypothesis in Theorem \ref{thm:essentiality_rational} cannot be removed: for
$M = \mathbb{Z}/p\mathbb{Z}$, every submodule localizes to zero, although $0$
is not essential in $M$.
\end{remark}

\section{Morita equivalence and factor rings}

We first consider the behavior of iso-multiplication modules under Morita
equivalence. Under the standard matrix equivalence, $M_A$ corresponds to a
row module, not to the regular module of the matrix ring.

\begin{theorem}
\label{thm:matrix_amplification} Let $A$ be an
arbitrary unital ring, $n \ge1$, $S = M_{n}(A)$, and $M_{A}$ any right module.
Endow $M^{1 \times n}$ with its natural right $S$-module structure by matrix
multiplication. Then
\[
M_{A} \text{ is im} \iff M_{S}^{1 \times n} \text{ is im}
\]
The same equivalence holds for classical multiplication modules.
\end{theorem}

\begin{proof}
Let $e = e_{11} \in S$. The standard Morita equivalence between right
$A$-modules and right $S$-modules sends $M$ to $F(M) = M^{1 \times n}$, with
inverse $X \mapsto Xe$; its exactness and its lattice bijection carry every
$A$-submodule $N \le M$ to the corresponding $S$-submodule $F(N) \le F(M)$.

Every two-sided ideal of $S$ has the form $M_{n}(B)$ for a unique two-sided
ideal $B$ of $A$. Direct calculation of the entries of row-matrix products
gives the equality of submodules
\[
F(M) M_{n}(B) = (MB)^{1 \times n} = F(MB).
\]
Suppose $M$ is iso-multiplication and take $X \le F(M)$. Write $X = F(N)$ under the
equivalence. For some $B \triangleleft A$, $N \cong MB$. Applying $F$ yields
$X \cong F(MB) = F(M)M_{n}(B)$; hence $F(M)$ is iso-multiplication. Conversely, if $F(M)$ is
iso-multiplication, apply the inverse equivalence to the isomorphism $F(N) \cong F(M)M_{n}(B)$
to obtain $N \cong MB$. Replacing isomorphisms by equalities and using the
submodule equality displayed above gives the identical argument for multiplication modules.
\end{proof}

\begin{corollary}
For any equivalence $F: \text{Mod-}A \to\text{Mod-}S$
of categories of right modules, $M_{A}$ is iso-multiplication if and only if $F(M)_{S}$ is iso-multiplication.
The same is true of the classical multiplication condition.
\end{corollary}

\begin{proof}
Choose a Morita progenerator ${}_{A} P_{S}$ realizing $F(M) \cong M
\otimes_{A} P$. For a two-sided ideal $B \triangleleft A$, let $J
\triangleleft S$ be the corresponding ideal under the Morita ideal-lattice
correspondence, characterized by $BP = PJ$. Since tensoring with $P$ is an exact Morita equivalence, the inclusion
$MB\hookrightarrow M$ induces an injection
$(MB)\otimes_A P\hookrightarrow M\otimes_A P$. Identifying its image
with a submodule of $M\otimes_A P$, we obtain
\[
(MB)\otimes_A P=(M\otimes_A P)J,
\]
Indeed, in $M\otimes_A P$ one has
$mb\otimes p=m\otimes bp$, so the image of $(MB)\otimes_A P$
is the submodule generated by the elements $m\otimes bp$ with
$b\in B$. Since $BP=PJ$, this is exactly the submodule generated
by $m\otimes pj$ with $p\in P$ and $j\in J$, namely
$(M\otimes_A P)J$. The equivalence also preserves and reflects
submodules and their isomorphisms. Thus the defining statement ``each submodule is isomorphic to an
ideal multiple'' is carried to precisely the corresponding statement over $S$.
Equalities yield the m-assertion.
\end{proof}

\begin{example}
Let $A = T_{2}(k)$ and $n \ge
1$, and set $S = M_{n}(A)$. By Example \ref{ex:basic_counterexample} and
Theorem \ref{thm:matrix_amplification}, the right $S$-module $A^{1 \times n}$
is iso-multiplication but is not a classical multiplication module. This is compatible with the fact that
the right regular module $S_{S}$ is not iso-multiplication for $n \ge2$: $A^{1 \times n}$ and
$S_{S}$ are different right $S$-modules. The latter failure follows from a
rank obstruction or, since $A$ is not semiprime, directly from the Morita
inverse $Se \cong A^{n}$ and the classification of modules over $A$: $A^{n}$
is not iso-multiplication when $n \ge2$.
\end{example}

\begin{corollary}
Let $S = M_{n}(T_{2}(k))$ and
let $F$ be the standard Morita equivalence $\text{Mod-}T_{2}(k) \to
\text{Mod-}S$. The iso-multiplication right $S$-modules, with no finiteness restriction, are
exactly the Morita images
\[
F(0), \quad F(S_{1}), \quad F(S_{2}), \quad F(P), \quad F(S_{1} \oplus S_{2}),
\quad F(P \oplus S_{2}).
\]
Among these six types, only $F(P \oplus S_{2})$ fails to be a classical
multiplication module.
\end{corollary}

\begin{proof}
An equivalence is essentially surjective and bijects isomorphism classes.
Combine Theorem \ref{thm:matrix_amplification} with Theorem
\ref{thm:T2_classification}, applying the last assertion of Theorem
\ref{thm:matrix_amplification} as well.
\end{proof}

\begin{remark}
Theorem \ref{thm:matrix_amplification} holds without finiteness or chain
conditions on the ring or the module. The proof uses the standard Morita
correspondence between submodules and two-sided ideals.
\end{remark}

\subsection{Factor rings and applications}

We give a factor-ring criterion and some consequences of the preceding results.

\paragraph{Passing through arbitrary factor rings.}

The following property is preserved when a module is viewed over a ring
and over a factor ring acting faithfully on it.

\begin{theorem}
\label{thm:factor_ring_inflation} Let $A$ be any unital ring, $H \triangleleft
A$ a two-sided ideal, $\overline{A} = A/H$, and $M$ a unitary right
$\overline{A}$-module, viewed also as an $A$-module. Then
\[
M_{A} \text{ is im} \iff M_{\overline{A}} \text{ is im}, \quad M_{A} \text{ is
an m-module} \iff M_{\overline{A}} \text{ is a multiplication module}.
\]
No commutativity, chain condition or finiteness assumption is needed.
\end{theorem}

\begin{proof}
Since $H$ annihilates $M$, the $A$-submodules and the $\overline{A}%
$-submodules of $M$ are identical. Similarly, an additive map between two such
submodules is $A$-linear if and only if it is $\overline{A}$-linear. For a
two-sided ideal $B$ of $A$, we have $MB = M(B+H)/H$, where the expression on
the right means multiplication by the two-sided ideal $(B+H)/H$ of
$\overline{A}$. Conversely, every two-sided ideal of $\overline{A}$ is the
image of a two-sided ideal of $A$ containing $H$. Thus the available ideal
multiples are exactly the same submodules, and their isomorphism types are
exactly those of the quotient action. The two equivalences follow separately
for equality and for abstract isomorphism.
\end{proof}

\begin{example}
Let $A = T_{n}(k)$, $n
\ge3$, and define the surjective homomorphism $\pi: A \to T_{2}(k)$ by taking
the upper-left $2 \times2$ block. Indeed the upper-left $2 \times2$ block of a
product of upper triangular matrices depends only on their upper-left blocks.
Take $M = T_{2}(k)_{T_{2}(k)}$, with the $A$-action $x \cdot a = x \pi(a)$. By
Theorem \ref{thm:factor_ring_inflation} and Example
\ref{ex:basic_counterexample}, $M_{A}$ is iso-multiplication but not a multiplication module. A witness is
the right submodule $ke_{22} \le T_{2}(k)$; it is isomorphic to $ke_{12}$, but
is not itself an ideal multiple. This construction proves existence for all
$T_{n}(k)$; it does not classify the iso-multiplication modules over $T_{n}(k)$.
\end{example}

\begin{example}
Let $0 \neq V$ be any $k$-vector
space, including $\dim_{k} V = 2$, and let $A(V)$ be the triangular algebra of
Theorem \ref{thm:triangular_family}. Its right regular module is iso-multiplication but not a multiplication module.
When $\dim_{k} V = 2$ this is the two-arrow (Kronecker) triangular algebra,
giving a direct example in a representation-infinite setting. The result is an
example for the regular module, not a classification of arbitrary Kronecker representations.
\end{example}

\paragraph{A sharpened consequence of Goldie rigidity.}

\begin{corollary}
\label{cor:reducedness_exact_rank} Suppose $A$ is semiprime right Goldie and
$A_{A}$ is iso-multiplication. Then
\[
Q_{\text{cl}}^{r}(A) \cong D_{1} \times\dots\times D_{s}, \quad A
\hookrightarrow D_{1} \times\dots\times D_{s}, \quad\text{u.dim}(A_{A}) = s,
\]
for division rings $D_{i}$. In particular, $A$ has no nonzero nilpotent
elements. If $A$ is prime, it is a domain and has right uniform dimension one.
\end{corollary}

\begin{proof}
Corollary \ref{cor:division_ring_factors} gives the displayed decomposition of
the classical right quotient ring. The classical localization map of a
semiprime right Goldie ring is injective, so $A$ embeds in this direct
product. Every nilpotent element in a division ring is zero, and the same
holds coordinatewise in a finite direct product, proving reducedness of $A$.
The right uniform dimension of $A$ is the length of the right regular
$Q$-module, by the standard Goldie rank identity. Each division-ring factor
contributes exactly one simple summand, so this dimension equals $s$. If $A$
is prime, its classical quotient ring is simple Artinian, hence there is only
one factor; its embedding in a division ring makes $A$ a domain.
\end{proof}

\begin{remark}
The conclusion that $Q$ is a finite product of
division rings is an obstruction to the iso-multiplication property, not a converse. Even when
$Q$ is a field, an $A$-module can have simple nonzero rationalization and
still fail iso-multiplication: the module $k[x] \oplus k[x]/(x)$ described in Example
\ref{ex:rational_not_sufficient} has precisely this feature. Its torsion
summand is invisible to rationalization but visible among integral submodule types.
\end{remark}

\paragraph{A direct length-spectrum falsification test.}

\begin{proposition}
\label{prop:missing_length_test} Let $M_{A}$ have
finite composition length $r$, and define
\[
\Lambda(M) = \{l_{A}(MB) : B \triangleleft A\}.
\]
If $M$ is iso-multiplication, then $\{0, 1, \dots, r\} \subseteq\Lambda(M)$. Consequently,
failure of even one integer $j \in\{0, \dots, r\}$ to belong to $\Lambda(M)$
proves that $M$ is not iso-multiplication.
\end{proposition}

\begin{proof}
Choose a composition series $0 = M_{0} \subsetneq M_{1} \subsetneq
\dots\subsetneq M_{r} = M$. The defining property gives an ideal $B_{i}
\triangleleft A$ with $M_{i} \cong_{A} MB_{i}$, and therefore $l_{A}(MB_{i}) =
i$ for each $i = 0, \dots, r$. Taking the contrapositive proves the test.
\end{proof}

\begin{example}
Let $A = k[\epsilon
]/(\epsilon^{2})$ and $M = k \oplus k$, with $\epsilon M = 0$. Then $l_{A}(M)
= 2$, while its ideal multiples are $M0=0$, $M(\epsilon)=0$, and
$MA=M$: the ideal $(\epsilon)$ annihilates $M$, whereas the unit ideal
$A$ acts as the identity. Thus $\Lambda(M) = \{0, 2\}$ and length one is missing. Hence $M$ is not
iso-multiplication. This calculation also illustrates that one need not classify all
submodules before applying the test.
\end{example}

\begin{corollary}
Suppose $A$ has exactly $t$
two-sided ideals, $M_{A}$ is iso-multiplication and $l_{A}(M) = t-1 < \infty$. Then the $t$
numbers $l_{A}(MB)$, as $B$ runs through the distinct two-sided ideals, are
precisely $0, 1, \dots, t-1$, each occurring exactly once.
\end{corollary}

\begin{proof}
Proposition \ref{prop:missing_length_test} produces all $t$ different lengths
$0, \dots, t-1$ among the ideal multiples, whereas there are only $t$ choices
for $B$. The pigeonhole principle proves the asserted one-to-one assignment.
\end{proof}

\begin{remark}
The missing-length criterion is an
immediately applicable formulation of Proposition \ref{prop:length_spectrum};
the sharpness corollary is a consequence of the same count. Neither is
advertised as a new independent classification. Likewise Corollary
\ref{cor:reducedness_exact_rank} is a strengthened exposition of the existing
Goldie result, not a second unrelated localization theorem.

The term \emph{isonoetherian} refers to stabilization of isomorphism types
along ascending chains, whereas \emph{isoartinian} refers to such
stabilization along descending chains \cite{Facchini2016}. The families above
show that iso-multiplication, unlike the classical ascending and descending chain conditions,
is an isomorphism-spectrum representation condition. In particular the
infinite-dimensional $A(V)$ example is iso-multiplication and isoartinian but not
isonoetherian (Corollary \ref{cor:iso_chain_separation}); finite-length
iso-multiplication modules satisfy both classical chain conditions, independently of the iso-multiplication definition.
\end{remark}

\paragraph{A brief commutative comparison.}

\begin{remark}
For commutative rings the two-sided-ideal
requirement is automatic, but the distinction between equality and isomorphism
of submodules remains logically meaningful. For modules annihilated by a power
of the maximal ideal over a commutative local ring, the iso-multiplication
condition reduces to classical multiplication, as follows from the local
cyclicity theorem above. For example, if $(R, \mathfrak{m})$ is commutative
local and $\mathfrak{m}^{h} M = 0$, then any nonzero cyclic submodule $Rx$
with $\mathfrak{m}^{h-1}x \neq0$ must be isomorphic to $IM$ with $I
\nsubseteq\mathfrak{m}$, so $I = R$ and $M \cong Rx$ is cyclic. Conversely,
every cyclic module over a commutative ring is a multiplication module.

This observation is included only to delimit the scope of the present
investigation, not as an independent main result. Existing studies of
classical multiplication modules, including the locally cyclic projective and
one-dimensional domain cases \cite{Ali2013, Tolooei2020}, are relevant to
further comparisons. We do not claim either a general equivalence or a
separating example for arbitrary commutative rings. Our main results concern
genuinely noncommutative examples, classifications, and obstructions.
\end{remark}

\section{Chain conditions and further remarks}

Recall that a module is isonoetherian (respectively, isoartinian) if the
terms of each ascending (respectively, descending) chain are eventually
isomorphic \cite{Facchini2016}. The iso-multiplication condition alone does
not compare ideal multiples by inclusion.

\begin{theorem}
Assume $M_{A}$ is iso-multiplication, and that (i) only
finitely many isomorphism types $[MB]$ occur as $B$ ranges over the two-sided
ideals of $A$, and (ii) whenever two submodules of $M$ embed in one another,
they are isomorphic. Then $M$ is both isonoetherian and isoartinian. In
particular, hypothesis (i) holds if $A$ has only finitely many two-sided ideals.
\end{theorem}

\begin{proof}
The inclusion chain of submodules induces a monotone chain in the finite
preorder of embeddability on isomorphism types. Condition (ii) makes the
preorder antisymmetric, hence a finite partial order. Every ascending and
descending chain in a finite partial order is eventually constant, so the
submodules are eventually isomorphic.
\end{proof}

\begin{remark}
Finiteness of the ideal lattice alone is not substituted for the
Schr\"oder--Bernstein condition (ii): different isomorphism types of infinite
modules may embed in one another. The more general embedding criterion below
keeps these logically independent requirements visible.
\end{remark}

\begin{remark}
For the triangular algebra
$A(V)$ with $V$ infinite-dimensional, Corollary \ref{cor:iso_chain_separation}
gives an iso-multiplication module which is isoartinian but not isonoetherian. Thus the two
iso-chain conditions must be investigated separately even within the present class.
\end{remark}

\begin{theorem}
If $D$ is a domain in
which every right ideal is principal, then $D_{D}$ is iso-multiplication, isonoetherian and
isoartinian. Moreover $D_{D}$ is a multiplication module if and only if every right ideal
is two-sided.
\end{theorem}

\begin{proof}
Every nonzero right ideal is $dD \cong D_{D}$ by left multiplication by $d
\neq0$; the only remaining type is $0$. In any ascending or descending chain,
the possible transitions between the two types $0$ and $D$ occur at most once.
The last claim follows because the ideal multiples $DB = B$ are precisely the
two-sided ideals.
\end{proof}

\begin{lemma}
\label{lem:matrix_size_obstruction} Let $D$ be a
principal right ideal domain and $n \ge2$. The regular right module of $S =
M_{n}(D)$ is not iso-multiplication.
\end{lemma}

\begin{proof}
Every two-sided ideal of $S$ is $M_{n}(I)$ for a two-sided ideal $I$ of $D$.
If $I \neq0$, then $I = dD$ for some $d \neq0$, so left multiplication by the
diagonal matrix $d1_{n}$ gives an injective right $S$-linear map $S \to S$
with image $M_{n}(dD) = M_{n}(I)$. Hence every nonzero two-sided ideal of $S$
is isomorphic to $S_{S}$. But the nonzero right ideal $e_{11}S$ is not
isomorphic to $S_{S}$: under the standard Morita equivalence, they correspond
to $D_{D}$ and $D_{D}^{n}$, respectively. A principal right ideal domain is a
right Ore domain and has well-defined rank over its classical division ring of
fractions, so these latter modules have different ranks when $n \ge2$. Thus
$e_{11}S$ cannot be isomorphic to any ideal multiple of $S$.
\end{proof}

\begin{theorem}
\label{thm:rigidity_facchini_nazemian} Let $A$ be a semiprime right noetherian
right isoartinian ring. Then $A_{A} \text{ is im}$ if and only if $ A \cong D_{1} \times\dots\times D_{s}$ for principal right ideal domains $D_{1}, \dots, D_{s}$. In this case $A_{A}$
is both isoartinian and isonoetherian.
\end{theorem}

\begin{proof}
	Assume first that \(A_A\) is an iso-multiplication module.
	By \cite[Theorem 4.7(1)]{Facchini2016}, under the stated hypotheses \(A\) is
	isomorphic to a finite direct product $	A\cong
	\prod_{i=1}^{s} M_{n_i}(D_i),$
	where each \(D_i\) is a principal right ideal domain.	By Theorem~\ref{thm:finite_direct_products},
	the iso-multiplication property of the regular module passes to each
	factor \(M_{n_i}(D_i)\). 
	Lemma~\ref{lem:matrix_size_obstruction} shows that
	\(n_i=1\) for every \(i\). Consequently, $	A\cong D_1\times\cdots\times D_s.$
	Conversely, suppose that $	A=\prod_{i=1}^{s}D_i,$
	where every \(D_i\) is a principal right ideal domain.
	Every right ideal of \(A\) has the form $	N=\prod_{i=1}^{s}N_i,$
	where \(N_i\) is a right ideal of \(D_i\).
	Since \(D_i\) is a domain and every right ideal is
	principal, each \(N_i\) is either zero or isomorphic
	to \(D_i\) as a right module. Let $	I=\prod_{i=1}^{s} I_i,
	$ where $
	I_i=
	\begin{cases}
		D_i,&N_i\ne 0,\\
		0,&N_i=0.
	\end{cases}$ \\
	Then \(I\) is a two-sided ideal of \(A\), and $	N\cong I=AI.$
	Thus every right ideal of \(A_A\) is isomorphic to
	an ideal multiple of \(A\), so \(A_A\) is iso-multiplication. Finally, consider any ascending or descending chain of right ideals of \(A\). Associate to each term the
	subset of indices \(i\) for which its \(i\)-th component
	is nonzero. These subsets are monotone along the chain.
	Since there are only finitely many indices, they
	eventually stabilize. Once they stabilize, all terms
	in the chain are isomorphic to the same product of
	nonzero factors \(D_i\). Hence \(A_A\) is both
	isonoetherian and isoartinian.
	
\end{proof}

\begin{remark}
The right-noetherian
assumption in Theorem \ref{thm:rigidity_facchini_nazemian} is part of the
cited classification, not an assumption that may silently be removed.
Nazemian-Yazdonov \cite{Nazemian2026} constructed a prime right and left
isoartinian ring which is neither right nor left noetherian (preprint dated 16
October 2026). That construction does not itself establish whether their ring
is iso-multiplication; it does decisively invalidate any attempt to assume semiprime
isoartinian rings must be noetherian.
\end{remark}

\begin{corollary}
If $A$ is a simple Artinian ring, then $A_{A}$
is iso-multiplication if and only if $A$ is a division ring.
\end{corollary}

\begin{proof}
Write $A = M_{n}(D)$ by Artin-Wedderburn. For $n \ge2$, use Lemma
\ref{lem:matrix_size_obstruction}; for $n = 1$, the only right ideals are $0,
D$.
\end{proof}

\subsection{An embedding criterion}

For an arbitrary right module $M$ let $\mathcal{T}(M) = \{[MB] : B
\triangleleft A\}$ be the collection of isomorphism classes of ideal
multiples, where $[X]$ denotes the class of $X$. Define $[X] \le[Y]$ when $X$
embeds as a submodule of $Y$. This is in general a preorder, rather than a
partial order: mutual embeddings need not imply isomorphism for arbitrary modules.

\begin{theorem}
Suppose $M$ is iso-multiplication
and the following hold:

\begin{enumerate}
\item[(i)] whenever $X = MB$ and $Y = MC$ embed in one another, then $X \cong
Y$;

\item[(ii)] the resulting partially ordered set of isomorphism types
$\mathcal{T}(M)$ has both the ascending and descending chain conditions.
\end{enumerate}

Then $M$ is both isonoetherian and isoartinian.
\end{theorem}

\begin{proof}
	Let $	N_1\subseteq N_2\subseteq\cdots$
	be an ascending chain of submodules of \(M\).
	Since \(M\) is an iso-multiplication module, for every \(i\) there is
	a two-sided ideal \(B_i\triangleleft A\) such that $	N_i\cong MB_i.$
	The inclusion \(N_i\hookrightarrow N_{i+1}\), together
	with these isomorphisms, gives an embedding $	MB_i\hookrightarrow MB_{i+1}.$
	Consequently, $	[MB_1]\preceq[MB_2]\preceq\cdots.$
	By condition (ii), this ascending chain of isomorphism
	classes eventually stabilizes. Hence there is an integer
	\(n_0\) such that $	[MB_i]=[MB_j]
$ $\text{for all }i,j\geq n_0.$
	It follows that $	N_i\cong MB_i\cong MB_j\cong N_j
$ $\text{for all }i,j\geq n_0.$
	Thus \(M\) is isonoetherian. The descending case is analogous. If $	N_1\supseteq N_2\supseteq\cdots,$
	then the corresponding ideal multiples give $	[MB_1]\succeq[MB_2]\succeq\cdots.$
	The descending chain condition in (ii) implies eventual
	stabilization of these isomorphism classes. Therefore
	\(M\) is isoartinian.
\end{proof}

\begin{remark}
This criterion is structural but conditional. It does not assert that the iso-multiplication
condition alone imposes a chain condition. The mutual-embedding hypothesis
cannot be omitted just by calling $\le$ a partial order. The $A(V)$ examples
show why explicit controls on ideal multiples are required.
\end{remark}

\subsection{Examples and limitations}

\begin{table}[htbp]
\centering
\par%
\begin{tabular}{@{}p{0.57\textwidth}cc@{}}
\toprule $M_{A}$ & m-module & im-module\\
\midrule A simple right module & Yes & Yes\\
$T_{2}(k)_{T_{2}(k)}$ & No & Yes\\
$K[x;\sigma]_{K[x;\sigma]}$ from Example \ref{ex:skew_polynomial} & No & Yes\\
$(T_{2}(k) \times k^{r})_{T_{2}(k) \times k^{r}}$ & No & Yes\\
$P = e_{11}T_{2}(k)$ & Yes & Yes\\
$S_{1} \oplus S_{2}$ over $T_{2}(k)$ & Yes & Yes\\
$S_{2} \oplus S_{2}$ over $T_{2}(k)$ & No & No\\
$D_{D} \oplus D_{D}$ for a division ring & No & No\\
\bottomrule
\end{tabular}
\caption{Comparison of multiplication modules and iso-multiplication modules}%
\label{tab:module_comparison}%
\end{table}

Theorem \ref{thm:T2_classification} applies to all modules over $T_2(k)$,
but not to arbitrary Artinian rings. In general, the iso-multiplication
property passes neither to submodules (Example
\ref{ex:submodule_failure}) nor to direct sums (Example
\ref{ex:direct_sum_failure}). The identity $N=M(N:M)$ for multiplication
modules also fails for iso-multiplication modules (Example
\ref{ex:basic_counterexample}).

\subsection{Connections with previous results}

Tuganbaev's definition of multiplication modules over noncommutative rings
requires the equality $N = MB$ with $B$ a two-sided ideal \cite{Tuganbaev2003,
Tuganbaev2004}. The present paper replaces equality by an $A$-linear
isomorphism, and Example \ref{ex:basic_counterexample} shows that the
relaxation is strict. Tuganbaev's results on distributivity, localization and
flatness are used only with their original hypotheses \cite{Tuganbaev2003,
Tuganbaev2005, Tuganbaev2006}; our Goldie theorem is formulated for semiprime
right Goldie rings and gives a necessary rational obstruction rather than a
converse criterion.

The normal form for right modules over $T_{2}(k)$ is the familiar
representation theory of the one-arrow quiver $A_{2}$. The assertion specific
to this paper is the ideal-multiple recognition in Theorem
\ref{thm:T2_classification}; the quiver decomposition is not claimed as new.
The matrix-amplification result uses standard Morita correspondence
\cite{Lam1999}. Alsuraiheed-Bavula \cite{Alsuraiheed2021} study direct sums of
classical multiplication modules, whereas our examples concern the weaker
iso-multiplication property.

Facchini-Nazemian introduced iso-Noetherian and iso-Artinian chain conditions
\cite{Facchini2016, Facchini2017}. The implication in Theorem
\ref{thm:rigidity_facchini_nazemian} uses their published classification under
semiprime, right-Noetherian and right-isoartinian hypotheses, with
matrix-factor elimination as an additional deduction. The chain conditions are
not built into our definition. For commutative locally cyclic projective
modules, the work of Ali \cite{Ali2013} and other classical classifications
\cite{Tolooei2020} explain why we do not present an extended commutative
classification here. The results above use these standard decomposition and equivalence theorems
with the hypotheses stated in their original forms.

\section{Conclusion}

We have constructed iso-multiplication modules which are not multiplication
modules and classified all such modules over $T_2(k)$. We have also proved
that the classical localization of an iso-multiplication module over a
semiprime right Goldie ring is multiplicity-free. Further results concern
composition length, local rings, Morita equivalence, and chain conditions
up to isomorphism. The localization condition is not sufficient in general.

\end{document}